\documentclass[prodmode,acmtomacs]{acmsmall} 
 
\usepackage[ruled]{algorithm2e}
\usepackage{amsmath}
\usepackage{ulem}

\newcommand{\maxent}{q}
\newcommand{\distrFamily}{\mathcal{G}}
\SetAlFnt{\small}
\SetAlCapFnt{\small}
\SetAlCapNameFnt{\small}
\SetAlCapHSkip{0pt}
\IncMargin{-\parindent}

\usepackage{subfig} 

\newtheorem{exampl}{Example}

\acmVolume{0}
\acmNumber{0}
\acmArticle{0}
\acmYear{2014}
\acmMonth{0}

\begin{document}

\markboth{Andreychenko, Mikeev, Wolf}{Model Reconstruction for Moment-based Stochastic Chemical Kinetics}

\title{Model Reconstruction for Moment-based Stochastic Chemical Kinetics}
\author{Alexander Andreychenko
\affil{Saarland University}
Linar Mikeev
\affil{Saarland University}
Verena Wolf
\affil{Saarland University}}

\begin{abstract}
Based on the theory of stochastic chemical kinetics, the inherent randomness and stochasticity of biochemical reaction networks can be accurately 
described by discrete-state continuous-time Markov chains.   
The analysis of such processes is, however,
computationally expensive and sophisticated numerical methods are required.
Here, we propose an analysis framework in which we integrate 
  a number of moments of the process instead  of the state probabilities.
  This results in a very efficient simulation of the time evolution of the process.
  In order to regain the   state probabilities from the moment representation, 
  we combine the fast moment-based simulation with a maximum entropy  approach for the 
  reconstruction of the underlying probability distribution.
  We investigate the usefulness of this combined approach in the setting of  
  stochastic chemical kinetics and present   numerical results for three
  reaction networks showing its efficiency and accuracy.
  Besides a simple 
  dimerization
  system, we study a bistable   switch system and
  a multi-attractor network with complex dynamics.
\end{abstract}




\acmformat{Alexander Andreychenko, Linar Mikeev,
and Verena Wolf, 2014. Model Reconstruction for Moment-based Stochastic Chemical Kinetics.}

\begin{bottomstuff}
 This research has been partially funded by  the German Research
     Council (DFG) as part of the Cluster of Excellence
     on Multimodal Computing and Interaction at Saarland University and
     the Transregional Collaborative
     Research Center ``Automatic Verification and Analysis of Complex
     Systems'' (SFB/TR~14~AVACS).

Author's address: Computer Science Department,
Saarland University, 66123 Saarbruecken, Germany;
\end{bottomstuff}

\maketitle

\section{Introduction}
\label{sec:intro}
During the last two decades discrete stochastic models have become a very popular
description of   biochemical reactions that take place in living organisms. 
They provide an appropriate representation of the discrete molecular populations in 
the cell and  accurately mimic the inherent randomness and discreteness of   molecular
interactions~\cite{fedoroff-fontana,mcadams-arkin99,thattai-oudenaarden,elowitz-etal02}.

The theory of stochastic chemical kinetics gives a rigorously justified stochastic description 
in terms of discrete-state continuous-time Markov chains~\cite{gillespie77}. 
The dynamics of the chain is governed by the Chemical Master Equation (CME) which describes the time evolution of the state probabilities.
However, the CME can be solved analytically only in a very limited number of cases.
The main difficulty arising in the numerical solution of the CME
is the curse of dimensionality: each chemical species that is involved in a reaction 
 adds one dimension to the state space of the Markov chain since the state of the chain
 is given by a population vector counting the current number of molecules for each species.
Tight bounds for molecular counts are usually not known a priori and thus 
the size of the state space that has to be considered is extremely large or even infinite
rendering the direct numerical integration of the CME  infeasible.
Sophisticated truncation approaches have been developed~\cite{Munsky06,FAUIET,MLEJournal}
which work well as long as the average population sizes remain small. 
The main idea is to concentrate on those population vectors containing a significant 
amount of probability mass. In this way only very rare 
behavior
of the process is neglected. 
If the reaction network contains highly abundant  chemical species then
the underlying probability distribution of the process becomes very large, even
when insignificant parts are truncated. In such a case it is advantageous  to change the 
representation of the distribution. 
The idea of methods of moments or 
moment closure methods 
is to replace the distribution  of the Markov chain
by its  moments up to a certain finite order~\cite{Engblom,StumpfJournal}.   
It is possible to derive differential equations that can be used to
approximate the time evolution of the moments. For instance, if the distribution of the 
chain is similar to a multivariate normal distribution, one can obtain a very accurate 
approximation of the distribution by tracking the average molecule counts and their
variances and covariances over time. For systems exhibiting more complex
behavior
such as oscillations or multi-modality,   moments of higher order are necessary for an accurate 
description~\cite{StumpfJournal}. For instance, the  moments of order three typically describe the skewness 
of the distribution while moments of order four are known to be a measure of the width of the tails
of the distribution. 

Often one is   interested in the probability
of certain events or in likelihoods of observations of the process.
However, usually   prior information regarding the properties of the distribution (e.g. approximately normally distributed) is not given and in such a case
regaining the probability distribution from the moment description is non-trivial. In fact 
it turns out that this problem, known as the classical moment problem,
 has a long history in other application domains and only 
recently very efficient methods for the reconstruction of the distribution became available.   

Given a number of moments of a random variable, there is in general no 
unique solution for the corresponding distribution.
However it is possible to define
a sequence of distributions that converges to the true one
whenever the number of constraints approaches infinity~\cite{mnatsakanov_recovery_2009}.  
Conditions for the existence of a solution  
are well-elaborated (such as Krein's and Carlemann's conditions)
but they do not provide a direct algorithmic way to create the reconstruction.
Therefore, Pade approximation~\cite{mead1984maximum}
and
inverse Laplace transform~\cite{Chauveau1994186} have been considered
but turned out to work only in restricted cases and require a
large number of constraints.
Similar difficulties are encountered when 
 lower and upper bounds for the probability distribution are derived~\cite{gavriliadis_moment_2008,tari_unified_2005,Kaas198687}.
Kernel-based approximation methods have been proposed where  one restricts  to a particular class of distributions~\cite{gavriliadis_truncated_2012,mnatsakanov_recovery_2009,chen_song_2000}.
The numerically most stable methods are, however, based on the 
maximum entropy principle which has its roots 
in statistical mechanics and information theory. The idea is to
choose from all distributions that 
fulfill
the moment constraints the 
distribution that 
maximizes the entropy. 
The maximum entropy reconstruction is the least biased estimate 
that 
fulfills the moment constraints
and it makes no assumptions about the missing information.
No additional knowledge about the
shape of the distribution neither a large number
of moments is necessary.
For instance, if only the first moment (mean) is provided 
the result of applying the maximum entropy principle is 
exponential distribution. In case of two moments
(mean and variance) 
the reconstruction is given by normal distribution.
Additionally, if   experimental data (or simulation traces) is available,
  data-driven 
maximum entropy methods can be applied~\cite{wu_weighted_2009,Golan1996559}.
Recently, notable progress has been made in the
development of numerical methods for the moment constrained maximum entropy problem
~\cite{abramov2010multidimensional,bandyopadhyay2005maximum,mead1984maximum},
where the main effort is put to the transformation
of the problem in order to overcome the numerical difficulties
that arise during the optimization procedure.

In this paper we propose a combination of the maximum entropy reconstruction
and the moment closure approach for the solution of the CME.
We approximate the moments over time and for a particular time instant
of interest we reconstruct the underlying distribution with a moment constrained maximum entropy approach. We do not make any
further assumptions about the distribution and study the feasibility,
efficiency and accuracy of this combined approach. 
To the best of our knowledge, the maximum entropy approach has
not yet been applied to stochastic models of biochemical reaction networks
before.
We consider three example networks which are small enough such that 
we can  compare our results with a nearly exact solution obtained via a direct numerical
integration of the CME.  
The maximum entropy approach
has been applied   to   chemical reaction networks in~\cite{Smadbeck}.
However, they restrict to finite-state models 
where
the entropy maximization becomes much easier since the 
support of the distribution is bounded. Here, we allow for infinite state 
space and present two infinite-state case studies. 
In addition, Smadbeck and Kaznessis compare
their results with statistical estimates obtained from Monte-Carlo 
simulations while we compare to results obtained via a numerical 
solution for which the approximation error is known.
Two of the examples that we consider have multi-modal distributions which
makes the reconstruction harder. However, 
our findings show that the combination of moment closure and 
maximum entropy reconstruction is surprisingly accurate also for complex systems and
it is very efficient
in terms of running times. In particular, the reconstruction part is very fast 
so that
the main advantage of the moment closure method - the short running time - remains, even if
it is combined with the maximum entropy approach. 
Thus, it provides a very useful alternative to
other analysis methods such as Monte Carlo simulations of the CME.
In particular, most methods do not scale in the number of molecules while
the efficiency of the moment closure approach is independent of the population sizes.
Short running times are particularly important if parameters of the process have to be
adjusted or if experiments must be designed~\cite{Lygeros}, since for such problems the 
model has to be 
analyzed for many different parameter combinations.

The paper is further organized as follows: 
we introduce our model in Section~\ref{sec:stoch}
and shortly explain how the CME can be numerically integrated to obtain accurate results for small systems.
In Section~\ref{sec:momcl} we discuss how the moment closure approach is applied to
the CME and in Section~\ref{sec:maxent} we describe the details of the maximum entropy approach 
and how it can be used to reconstruct the distribution from a number of moments.
Finally, we present experimental results for the three case studies in Section~\ref{sec:exp} and 
conclude the paper with Section~\ref{sec:conc}.

\section{Stochastic Chemical Kinetics}\label{sec:stoch}
Stochastic chemical kinetics refers to a widely-used modelling framework for
the description of networks of biochemical reactions~\cite{mcquarrie67}.
We consider a biological compartment (e.g. a living cell) in which
molecules of different types undergo chemical reactions. Assuming that
this reaction volume is well-stirred and in thermal equilibrium, it is possible
to physically justify a Markov chain description of the chemical populations~\cite{gillespie77}, that is,
we consider a random vector $X(t)=(X_1(t),\ldots,X_n(t))$ where $X_i(t)$ is the
number of molecules of type $i$ at time $t$ ($i\in\{1,\ldots,n\}, t\ge 0$).
We assume that the set of possible reactions is given by the stoichiometric
equations
$$
R_j:\quad\quad  \ell_{j,1} S_1+\ldots+\ell_{j,n} S_n \xrightarrow{c_j}\tilde\ell_{j,1} S_1+\ldots+\tilde\ell_{j,n} S_n,
\quad 1\le j \le m.
$$
Here, $\ell_{j,i}$ and $\tilde\ell_{j,i}$ refer to the number of molecules used up and produced by the reaction,
respectively ($\ell_{j,i},\tilde\ell_{j,i}\in\mathbb N_0$ and $1\le i\le m$) and $c_j$ is the so-called stochastic
reaction rate constant that determines the probability of the reaction as explained in the sequel.

\begin{exampl}
\label{ex:simple_dimerization_reactions}
We consider a simple dimerization example~\cite{StumpfJournal} with the stoichiometric equations
\begin{alignat}{3}
R_1: &\quad2 P &\ \stackrel{c_1}{\longrightarrow}\ P_2,     \nonumber\\
R_2:  &\quad P_2 &\ \stackrel{c_2}{\longrightarrow}\ 2 P,     \nonumber\\
\end{alignat}
where  $\ell_{1,1}=\tilde\ell_{2,1}=2$, $\ell_{2,2}=\tilde \ell_{1,2}=1$ and $\ell_{1,2}=\ell_{2,1}=\tilde\ell_{1,1}=\tilde\ell_{2,2}=0$. Note
that we omit terms which are zero and the factors equal to one.
\end{exampl}

\subsection{Transition Rates}

Let $v_j\in \mathbb Z^n$ be the vector that describes the
population change of reaction $R_j$, that is, $v_j=(\tilde\ell_{j,1}-\ell_{j,1},\ldots,\tilde\ell_{j,n}-\ell_{j,n})$.
Transitions of the  Markov chain
$X$ correspond to chemical reactions and
the transition rate of reaction $R_j$ is given by
$$
\lim_{h\to 0} \frac{1}{h} P(X(t+h)=x+v_j\mid X(t)=x)=c_j \prod_{i=1}^n {x_i\choose \ell_{j,i}}
$$
where $x=(x_1,\ldots,x_n)\in\mathbb N_0^n$ is a state and  the binomial coefficients describe the number
of possible combinations of reactant molecules. Note that $c_j$ depends on the physical
properties of the reactants as well
 as on the temperature and the size of the reaction volume.
 Here we make the usual assumption that $c_j$ is constant in time.
 It is possible to extend the
results presented in the sequel for time-dependent $c_j$.
In this case the underlying Markov chain is
time-inhomogeneous and an accurate simulation may be
challenging if the rates strongly vary in time~\cite{qapl11}.

In the sequel we also restrict to chemical reactions that are at most bimolecular, i.e., we assume
that $\sum_{i=1}^n\ell_{j,i} \in \{0,1,2\}$, which is a reasonable assumption because reactions
where more than two molecules have to collide can usually be decomposed into smaller ones where
at most two molecules have to collide~\cite{gillespie77}.

\subsection{Chemical Master Equation}
The dynamics of $X$ is given by the chemical master equation that describes the time evolution of
the transient distribution $\pi(x,t)=P(X(t)=x)$ as a linear ordinary differential equation
\begin{equation}
\label{eq:master}
\frac{\partial \pi(x,t)}{\partial t}
= \sum_{j=1}^m \left(\alpha_j(x-v_j)\pi(x-v_j,t)
 -\alpha_j(x)\pi(x,t)\right).
 \end{equation}
Here, $\alpha_j(x)=c_j \prod_{i=1}^n {x_i\choose \ell_{j,i}}$ is the transition rate in state $x$ for reaction $R_j$.
 If an initial distribution, say at time $t=0$, is given then the equation
has a unique solution at all  finite times $t \ge 0$.
 It is important to note that often the number of states with $\pi(x,t)>0$ is infinite since
 bounds on the population sizes are not known a priori.  Thus, although in reality the molecule numbers are
 always finite, theoretically an infinite number of states can have positive probability.
 This leads to two complications compared to finite-state models.
First, the limiting distribution of the Markov chain may not exist and additional conditions are necessary
to ensure the existence~\cite{NLAA}. And second, truncation techniques are necessary to numerically simulate \eqref{eq:master} since only a tractable number of states can be considered in each integration step.
Only in special cases an analytic solution of \eqref{eq:master} is available~\cite{Jahnke_CME_analytical}.

\subsection{Direct Numerical Simulation}\label{sec:directnumsim}
In the sequel we shortly explain how the master equation in \eqref{eq:master} can be simulated numerically
since our goal is to compare a moment closure approximation combined with a reconstruction of the
individual probabilities with such a direct numerical simulation. The latter 
performs well as long as
the average population sizes remain small and the approximation error can be controlled by a simple
threshold criterion. Thus, we will be able to determine the accuracy of the moment closure approximation
as well as the accuracy of the reconstruction algorithm.

The direct numerical simulation that we consider is based on the dynamic state space truncation
developed for uniformization methods~\cite{IETFAU} and for integration schemes
such as Runge Kutta methods~\cite{MLEJournal,valuetools11}.
The main idea is to exploit the inflow-outflow form of \eqref{eq:master} for the construction of
the dynamic state space. The terms $ \alpha_j(x-v_j)\pi(x-v_j,t)$ can be seen as the inflow
to state $x$ for reaction $R_j$ while $\alpha_j(x,t)\pi(x,t)$ is the corresponding outflow.
Let $p(x,t)$ be the approximation of $\pi(x,t)$ during the numerical integration for all $x$
and all $t\ge 0$.
Initially we set $p(x,0)=\pi(x,0)$ and during an 
integration step for the interval $[t,t+h)$ we start  with a subset $S^{(t)}$ of states
 that have significant probability at time $t$, i.e.,
 $$S^{(t)}:=\{x\mid p(x,t)>\delta_1\}$$ where $\delta_1>0$ is a small threshold.
 For all states not in  $S^{(t)}$ we let $p(x,t)=0$.
 During the numerical integration we add new states to $S^{(t)}$ whenever they
 receive a significant amount of inflow, i.e. if we use the explicit Euler method,
 the new state probability at time $t+h$ is calculated as
 $$
 p(x,t+h)=p(x,t)+h\cdot \sum_{j=1}^m \left(\alpha_j(x-v_j)p(x-v_j,t)
 -\alpha_j(x)p(x,t)\right).
 $$
 For a state $x\not\in S^{(t)}$ this reduces to
  $$
 p(x,t+h)= h\cdot \sum_{j=1}^m  \alpha_j(x-v_j)p(x-v_j,t).
 $$
 Hence, we can loop over all states in $S^{(t)}$ and, before integrating their
 probability,  check whether their successors receive significant inflow. More
 precisely, we simply add a state $x$ to the set $S^{(t)}$ if
 $h\cdot \alpha_j(x-v_j)p(x-v_j,t)>\delta_2$ for some $j$. 
 Here, $\delta_2$ is again a small threshold.
 We then also compute $ p(x,t+h)$ for this new state.
It turns out that for most example networks  an accurate approximation is obtained if
we   work with a single threshold $\delta_1=\delta_2=:\delta$ and
 choose $\delta \in\{10^{-10},10^{-9}\ldots,10^{-5}\}$.
 Note that the new set $S^{(t+h)}$ will then contain all states $x\in S^{(t)}$
 whose probability at time $t+h$ is at least $\delta$ as well as all successors
 $x+v_j \not\in S^{(t)}$  where $x\in S^{(t)}$ and there exists a $j$ such that
 $h\cdot \alpha_j(x-v_j)p(x-v_j,t)>\delta$ (which implies that their probability at time
 $t+h$ is at least $\delta$). Obviously, different truncation strategies are possible
 (e.g. choose $\delta_2$ smaller than $\delta_1$). However, we found that simply adding
 all successors ($\delta_2=0$) is not efficient since often we have reversible reactions, i.e.,
 $v_j=-v_k$ for some $j\neq k$ where one direction is much more likely then the other, say $R_j$.
 In such a case the main part of the probability mass moves in the direction of $v_j$
 and the accuracy gain in adding a successor w.r.t. $v_k$ is not worth the effort since during
 the next construction of $ S^{(t)}$ these successors are anyway removed from $ S^{(t)}$.

In order to illustrate the method, we list  the size of the truncated state space and
the total loss of probability mass for the
  following example:

\begin{exampl}
\label{ex:exSwitchReactions}
	We consider a gene regulatory network called the exclusive switch \cite{loinger-lipshtat-balaban-biham07}.
	It describes the dynamics of two genes with an overlapping promotor region,
	and their products $P_1$ and $P_2$.
	Molecules of both species $P_1$ and $P_2$ are produced if no transcription factor is bound to the promotor region (region is free).
	However if a molecule of type $P_1$ $(P_2)$ is bound to the promotor then
	it inhibits the expression of the other product,
	i.e. only molecules of $P_2$ $(P_1)$ can be produced. Only one molecule
	can be bound to the promotor region at a time. 
	The model has an infinite state space and the stoichiometric equations are given by:
	\begin{equation*}
	\begin{aligned}[c]
		DNA &\ \stackrel{c_1}{\longrightarrow}\ DNA + P_1\\
		DNA &\ \stackrel{c_2}{\longrightarrow}\ DNA + P_2\\
		P_1 &\ \stackrel{c_3}{\longrightarrow}\ \emptyset \\
		P_2 &\ \stackrel{c_4}{\longrightarrow}\ \emptyset \\
	\end{aligned}
	\qquad
	\begin{aligned}[c]
		DNA + P_1 &\ \stackrel{c_5}{\longrightarrow}\ DNA.P_1\\
		DNA + P_2 &\ \stackrel{c_6}{\longrightarrow}\ DNA.P_2\\
		DNA.P_1 &\ \stackrel{c_7}{\longrightarrow}\ DNA + P_1\\
		DNA.P_2 &\ \stackrel{c_8}{\longrightarrow}\ DNA + P_2\\
	\end{aligned}
	\qquad
	\begin{aligned}[c]
		DNA.P_1 &\ \stackrel{c_9}{\longrightarrow}\ DNA.P_1 + P_1\\
		DNA.P_2 &\ \stackrel{c_{10}}{\longrightarrow}\ DNA.P_2 + P_2\\
	\end{aligned}
	\end{equation*}
	where the reaction rate constants $c_1,\ldots, c_{10}$ are given by the entries of the vector $c=(2.0,5.0,0.005,0.005,0.005,0.002,0.02,0.02,2.0,5.0)$
	and
	the initial conditions are such that only one $DNA$ molecule is present in the system 
	while the molecular counts for the rest of species are zero.

\begin{figure}[t]
\centering
  \subfloat[][]{\includegraphics[width=0.33\textwidth]{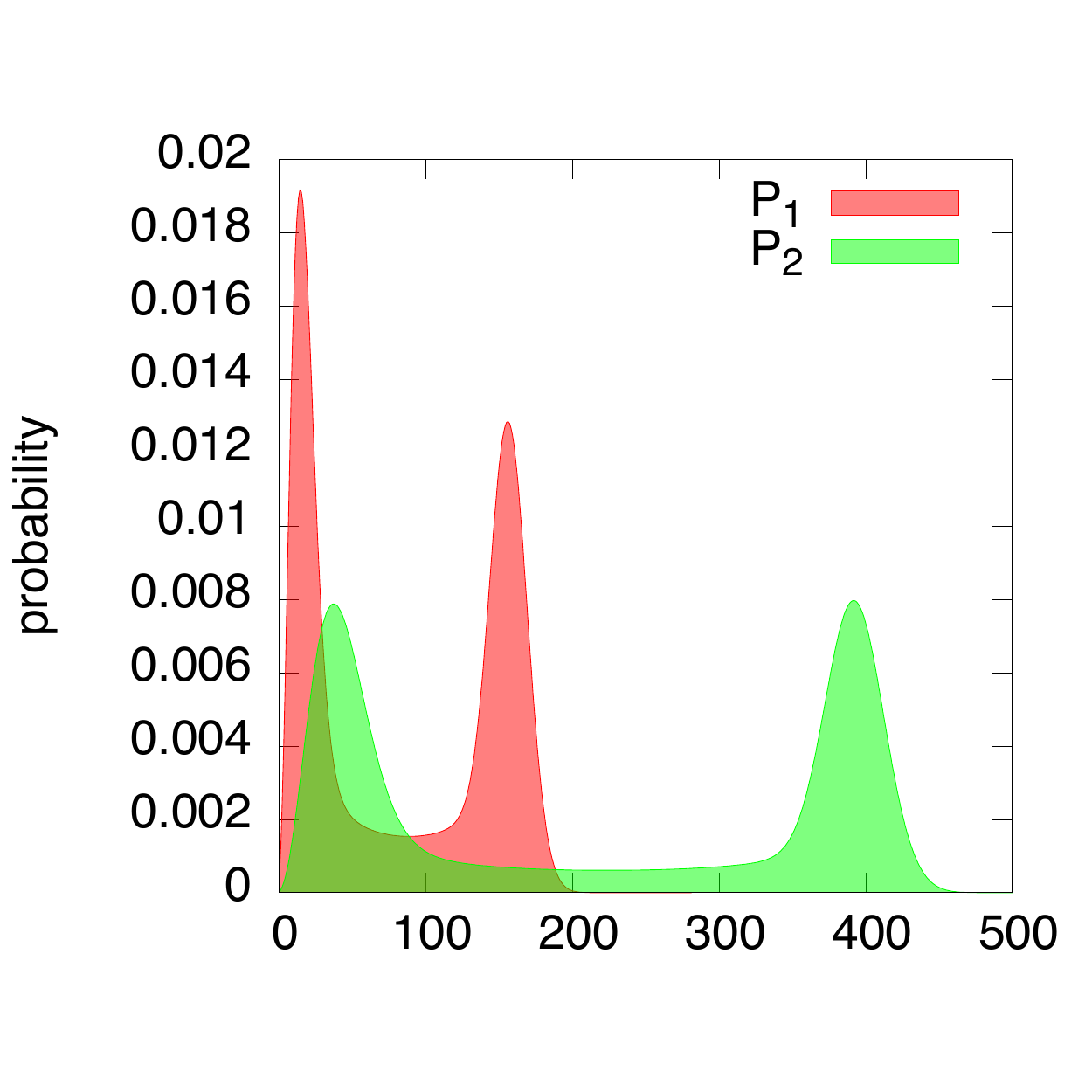}}
  \subfloat[][]{\includegraphics[width=0.33\textwidth]{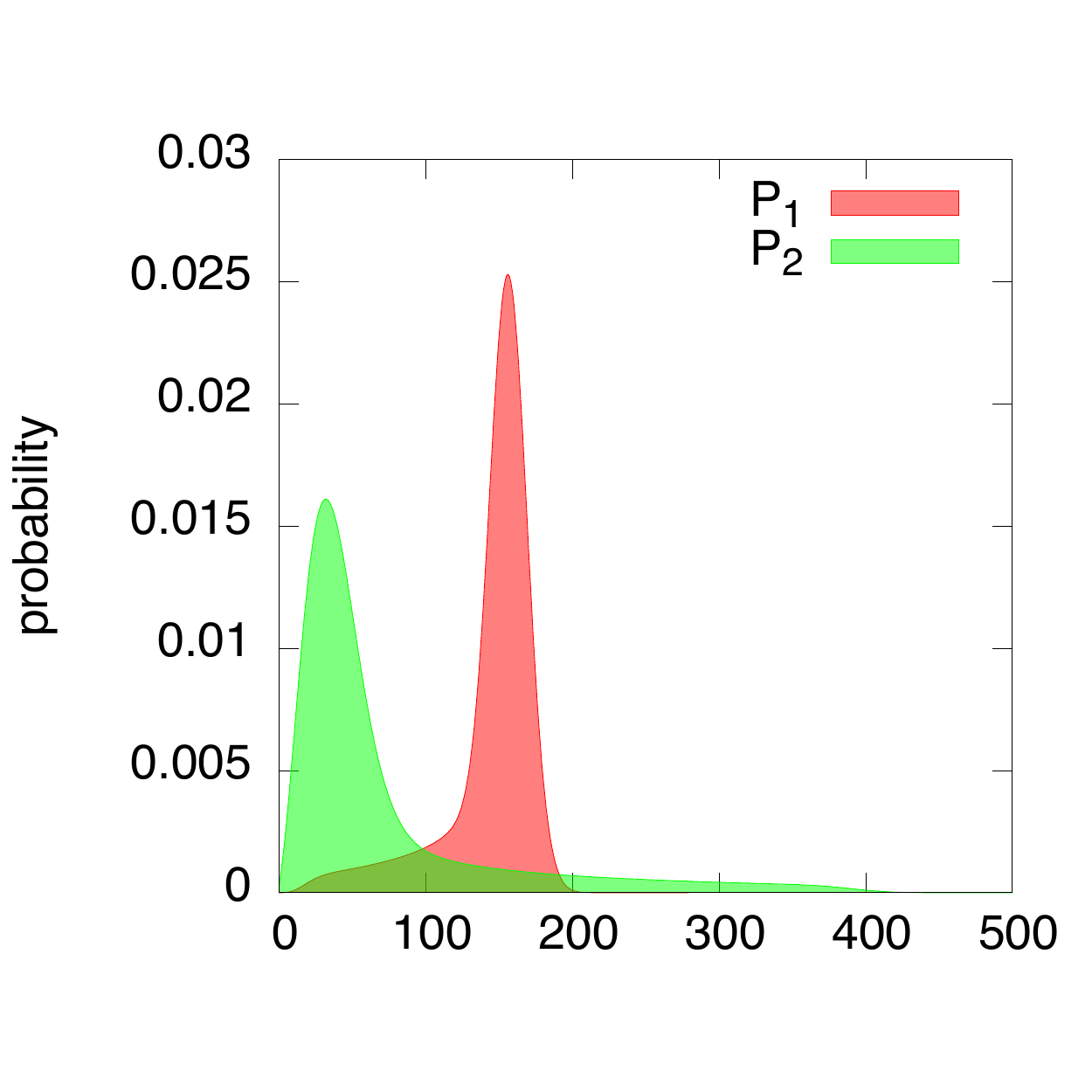}}
  \subfloat[][]{\includegraphics[width=0.33\textwidth]{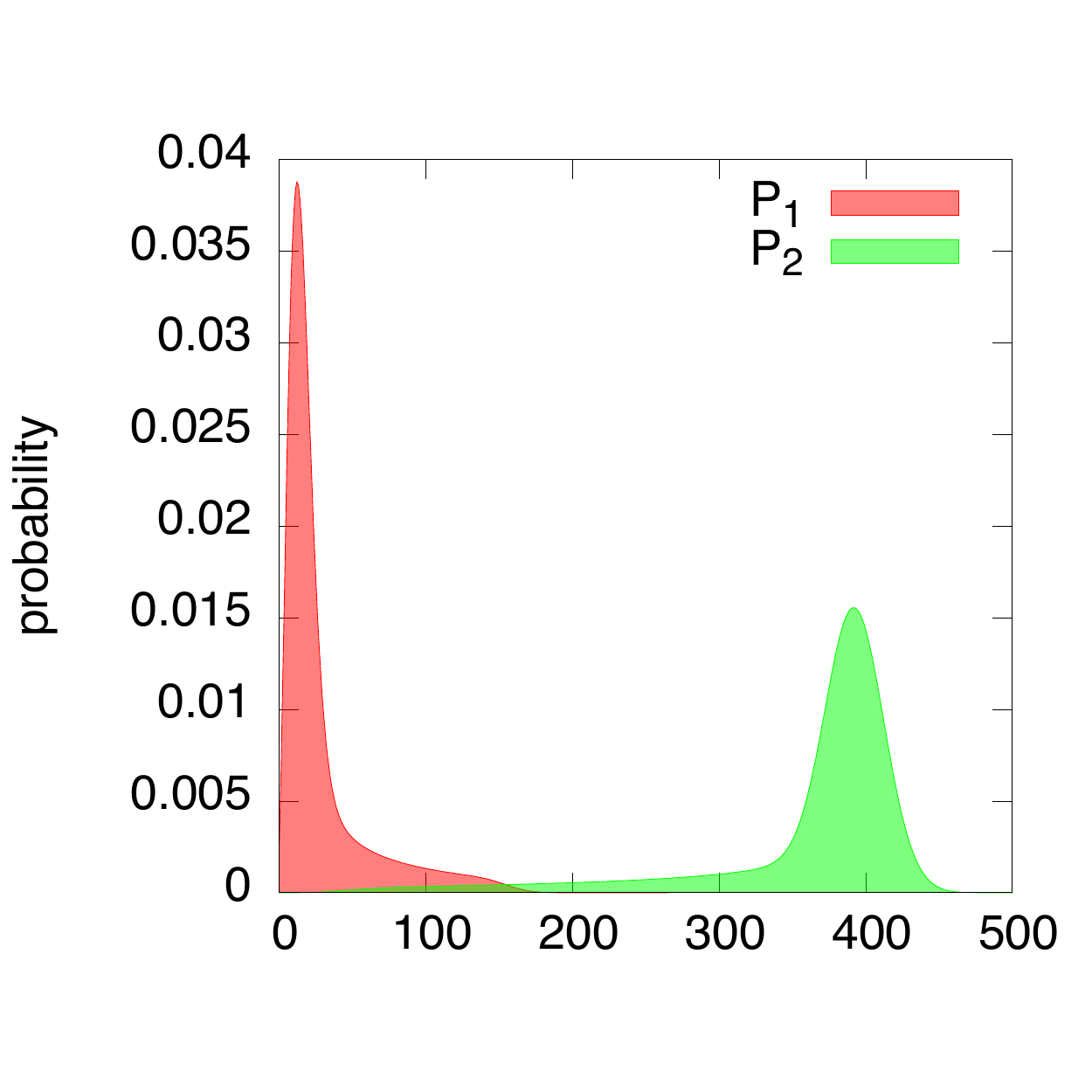}}
	\caption{\label{fig:exswitch_distr}Probability distribution of the protein counts
	$P_1$ and $P_2$ conditioned on the events that the promotor region is (a)  free (b) bound to $P_1$, and (c) bound to $P_2$,
	computed at time instant $t=100$ for exclusive switch system.}
\end{figure}

In Figure~\ref{fig:exswitch_distr} we plot the results of a direct numerical simulation
using the dynamical state space as explained above.
The different subfigures show
the marginal distributions of protein counts $P_1$ and $P_2$
when we condition on the three different states of the promotor region (free, $P_1$ or $P_2$ bound).
\renewcommand{\arraystretch}{1.2}
\begin{table}[b]
\tbl{Dynamical state space truncation results for the exclusive switch\label{tab:exswitch_truncation_error}}{
\begin{tabular}{|c|c|c|c|}
\hline 
$\delta$ & $|S|$ & $\epsilon$ & time (sec)\\
\hline
$10^{-10}$ & 183210 & $3\cdot 10^{-6}$ & 154\\
\hline
$10^{-12}$ & 203948 & $2\cdot 10^{-8}$ & 174\\
\hline
$10^{-15}$ & 265497 & $9\cdot 10^{-11}$ & 239\\
\hline
$10^{-20}$ & 381374 & $1\cdot 10^{-13}$ & 1027\\
\hline
\end{tabular}}
\end{table}
To investigate the accuracy of the obtained results
we refer to the Table~\ref{tab:exswitch_truncation_error}, where we list
the amount $\epsilon$ of   probability mass  
lost during the computation and 
the 
size $|S|$ of the truncated state space $S^{(t)}$ for different thresholds $\delta$
at time instant $t=100$.
Note that the probability of all states not in  $S^{(t)}$ is approximated
as zero. Thus, $\epsilon$ is equal to the total approximation error (sum of all
state-wise errors)
and all computed state probabilities are underapproximations,
i.e. \begin{equation}\label{eq:eps}
\sum_x  \pi(x, t) - \sum_x  p(x, t)\le \epsilon \end{equation}
where we have equality at the final time instant of the computation.
We find that, for instance, if we choose $\delta=10^{-15}$ the total approximation error
remains below $10^{-10}$.  \end{exampl}



\section{Moment Closure Approximation}\label{sec:momcl}
\newcommand{\Ex}[1]{\ensuremath{\textstyle E\left(#1\right)}}
As opposed to the method in the previous section, during the moment closure approximation
we     integrate the first $k$ moments of the distribution $\pi(x,t)$ over time.
For this we derive a system of differential equations for the moments along
the lines of Ale et al. to show how and where approximations
errors occur~\cite{StumpfJournal}. We restrict ourselves to the first two moments
in order to keep this review of moment closure techniques short and also because the
derivation of the equations for the first two moments are sufficient for illustrating
the technique.

 Let $f:\mathbb N^n_0\to \mathbb R^n$ be a function that is independent of $t$.
 In the sequel we will exploit the relationship
\begin{equation}\label{eq:expgen}
\begin{array}{lcl}
  \frac{d}{dt}\Ex{f(X(t))} &=& \sum\limits_{x}
f(x)\cdot  \frac{d}{dt} P(X{(t)}=x)\\[2ex]
&=& \sum\limits_{j=1}^m \Ex{\alpha_j(X(t)) \cdot (f(X(t) +  v_j ) - f(X(t)))}.
\end{array}
\end{equation}
 For $f(x)=x$ this yields a system of equations for the population means
\begin{equation}\label{eq:mean1}
\begin{array}{lcl}
  \frac{d}{dt}\Ex{X(t)} &=&   \sum\limits_{j=1}^m v_j \Ex{\alpha_j(X(t))}.
\end{array}
\end{equation}
Note that for bimolecular reactions, $\Ex{\alpha_j(X(t))}$ is equal to
$$ \begin{array}{ll}
c_j \Ex{X_{i}(t)\cdot X_{i'}(t)} \neq c_j \Ex{X_i(t)}\cdot \Ex{X_{i'}(t)} & \mbox{ if }i\neq i',\\[1ex]
0.5 c_j \Ex{X_{i}(t)\cdot (X_{i'}(t)-1)} \neq 0.5 c_j \Ex{X_i(t)}\cdot (\Ex{X_{i'}(t)}-1) & \mbox{ if }i=i',
\end{array} $$ for  
$i,i' \in \{1,\ldots,n\}$, 
which means that the system of ODEs in Eq.~\eqref{eq:mean1} is only closed
if at most monomolecular reactions ($\sum_{i=1}^n\ell_{j,i}\le 1$) are involved.
For most networks  the latter condition is not true
but we can approximate the unknown term $\Ex{X_i(t)\cdot X_{i'}(t)}$
either by assuming that the covariance is zero, which gives
$\Ex{X_i(t)\cdot X_{i'}(t)}=\Ex{X_i(t)}\cdot \Ex{X_{i'}(t)}$ or by
extending the system in~\eqref{eq:mean1}  with additional
equations for the second moments.
The general strategy is to replace $\alpha_j(X(t))$ by a Taylor series about
the mean $\Ex{X(t)}$. Let us write $\mu_i(t)$ for $\Ex{X_i(t)}$ and $\mu(t)$
for the vector with entries $\mu_i(t)$, $1\le i\le n$. Then
 $$
 \begin{array}{rcl}
 \Ex{\alpha_j(X)}& =& \alpha_j(\mu) + \frac{1}{1!}  \sum_{i=1}^n\Ex{X_i-\mu_i}\frac{\partial}{\partial x_i}\alpha_j(\mu)\\[1ex]
 &&+ \frac{1}{2!} \sum_{i=1}^n \sum_{k=1}^n\Ex{(X_i-\mu_i)(X_k-\mu_k)} \frac{\partial^2}{\partial x_i \partial x_k}\alpha_j(\mu)+\ldots
 \end{array}
 $$
 where we omitted $t$ in the equation to improve the readability.
 Note that $\Ex{X_i(t)-\mu_i}=0$ and since we restrict to reactions that are at most bimolecular, all terms of order three and more disappear. By letting $C_{ik}$ be the covariance
 $\Ex{(X_i(t)-\mu_i)(X_k(t)-\mu_k)}$ we get
 \begin{equation}\label{eq:meanrate}
 \begin{array}{rcl}
 \Ex{\alpha_j(X)}& =& \alpha_j(\mu) +  \frac{1}{2} \sum_{i=1}^n \sum_{k=1}^n C_{ik}\frac{\partial^2}{\partial x_i \partial x_k}\alpha_j(\mu)
 \end{array}
 \end{equation}
 Next, we derive an equation for the covariances by first exploiting the relationship
\begin{equation}\label{eq:cov}
 \frac{d}{dt} C_{ik} =  \frac{d}{dt} \Ex{X_iX_k}- \frac{d}{dt}(\mu_i\mu_k) =
  \frac{d}{dt} \Ex{X_iX_k}- (\frac{d}{dt}\mu_i)\mu_k-\mu_i(\frac{d}{dt}\mu_k)
 \end{equation}
 and if we couple this equation with the equations for the means, the only unknown term
 that remains is the derivative $\frac{d}{dt} \Ex{X_iX_k}$ of the second moment.
 For this we can use the same strategy as before, i.e., from Eq.~\eqref{eq:expgen}
 we get
  \begin{equation}\label{eq:secondm}
 \frac{d}{dt} \Ex{X_iX_k}  =    \sum_{j=1}^m \left[v_{j,i}v_{j,k}  \Ex{\alpha_j(X)}
 + v_{j,k}\Ex{\alpha_j(X)X_i}+ v_{j,i}\Ex{\alpha_j(X)X_k}\right]
 \end{equation}
 where $v_{j,i}$ and $v_{j,k}$ are the corresponding entries of the vector $v_j$.
 Clearly, we can use Eq.~\eqref{eq:meanrate} for the term $\Ex{\alpha_j(X)}$ while
   the terms $\Ex{\alpha_j(X)X_i}$ and $\Ex{\alpha_j(X)X_k}$ have to be
   replaced by the corresponding Taylor series about the mean.
  Let $f_j(x):=\alpha_j(x)x_i$. Similar to Eq.~\eqref{eq:meanrate} we get
  \begin{equation}\label{eq:meanrate2}
   \begin{array}{rcl}
 \Ex{\alpha_j(X)X_i}& =& \alpha_j(\mu) \mu_i+ \frac{1}{1!}  \sum_{i=1}^n\Ex{X_i-\mu_i}\frac{\partial}{\partial x_i}f_j(\mu)\\[1ex]
&&+ \frac{1}{2!} \sum_{i=1}^n \sum_{k=1}^n\Ex{(X_i-\mu_i)(X_k-\mu_k)} \frac{\partial^2}{\partial x_i \partial x_k}f_j(\mu)+\ldots.
 \end{array}
 \end{equation}
  Here, it is important to note that moments of order three come into play since derivatives
 of order three of $f_j(x)=\alpha_j(x)x_i$ may be nonzero. It is possible to take these terms into account
 by deriving additional equations for moments of order three and higher.
 Obviously, these equations will then include moments of even higher order such that
 theoretically we end up with an infinite system of equations. However, a popular
 strategy is to close the equations by assuming that all moments of order  $>M$
  that are centred around the mean are equal to zero. E.g. if we choose $M=2$, then
 we can simply use the approximation
   $$
   \begin{array}{rcl}
 \Ex{\alpha_j(X)X_i}& \approx & \alpha_j(\mu) \mu_i
 + \frac{1}{2!} \sum_{i=1}^n \sum_{k=1}^n\Ex{(X_i-\mu_i)(X_k-\mu_k)} \frac{\partial^2}{\partial x_i \partial x_k}f_j(\mu).
 \end{array}
 $$
 This approximation is then inserted into Eq.~\eqref{eq:secondm}
 and the result replaces  the term  $\frac{d}{dt} \Ex{X_iX_k}$ in Eq.~\eqref{eq:cov}.
 Finally, we can integrate the time evolution of the means and that of the covariances
 and variances.

\begin{exampl}
\newcommand{\fst}{\ensuremath{1}}
\newcommand{\snd}{\ensuremath{2}}
To illustrate the method we consider again the simple 
dimerization 
reaction 
system of Example~\ref{ex:simple_dimerization_reactions}.
Assuming that all central moments of order three and higher are equal to
zero, we get the following equations for the means, variances
and the covariance of the species
$$
\begin{array}{ll}
	\frac{d}{dt} \mu_{\fst}&=
	  -c_1  \mu_{\fst} \left( \mu_{\fst}-1 \right)
		-c_1 C_{\fst,\snd}
		+ 2 c_2 \mu_{\snd}\\[1ex]
	\frac{d}{dt} \mu_{\snd}&=
		  \frac{1}{2} c_1 \mu_{\fst} \left( \mu_{\fst}-1 \right)
		 + \frac{1}{2} c_1 C_{\fst,\snd}
		 - c_2 \mu_{\snd} \\[1ex]
  \frac{d}{dt} C_{\fst,\fst} &=
  -2 c_1 \mu^3_{\fst} + 4 c_1 \mu_{\fst}^2
	-2 c_1 \mu_{\fst} + 4 c_2 \mu_{\snd} + 4 c_2 \mu_{\fst} \mu_{\snd} \\[1ex]
	\frac{d}{dt} C_{\fst,\snd} &=
  	-\frac{3}{2} c_1 C_{\fst,\fst} 
	   +\mu_{\fst}\left(-c_1(\mu_{\fst}-1)-\frac{1}{2}c_1\mu_{\fst}
	   C_{\fst,\fst}+c_2\mu_{\fst}^2\right)\\
	   &+\mu_{\snd} \left(-2c_2-c_2\mu_{\fst}+c_1 C_{\fst,\fst}\right)\\[1ex]
	\frac{d}{dt}C_{\snd,\snd} &=
		\frac{3}{2} c_1 C_{\fst,\fst}
		+ \frac{1}{2}c_1 \mu_{\fst} (\mu_{\fst}-1)
		+ \mu_{\snd} (c_2 - c_1 C_{\fst,\fst})
\end{array}
$$
where we denote the expectations of species $P$ $(P_2)$ by 
$\mu_{\fst}$ ($\mu_{\snd}$),
variances are given by $C_{\fst,\fst}$, $C_{\snd,\snd}$ and 
the covariance between $P$ and $P_2$ is $C_{\fst,\snd}$.
In the equations we omit   $t$ to improve readability.

\end{exampl}

 In Section~\ref{sec:exp} we study the accuracy of the above example (cf. Table~\ref{tab:dimerization_mc_error})
and find that the approximation provided by moment closure method is very accurate even if only 
the means and covariances are considered. 
In general, however, experimental results show that the approximation tends to become worse if
systems that exhibit complex behavior such as multistability or oscillations.
Increasing the number of moments typically improves the accuracy~\cite{StumpfJournal}
but sometimes the resulting equations may become very stiff~\cite{Engblom}.

Grima has investigated the accuracy of the approximation for $n=2$ and $n=3$
by a comparison with the system size expansion of the master equation~\cite{grima2012study}.
He found that for monostable systems with large volumes  the approximation of the means
$\mu(t)$ have a relative error that scale as $\Omega^{-n}$
while the relative errors of the variances  and covariances scale as $\Omega^{-(n-1)}$, $n\in\{2,3\}$.
For small volumes or systems with multiple modes, however, only experimental evaluations
of the accuracy are available~\cite{StumpfJournal,Engblom}, where the approximated moments are compared to
statistical estimates based on Monte Carlo simulations of the process.
In Section~\ref{sec:exp} we focus on experimental results for the reconstructed probability distribution.
 However,  we also  compare the moments approximated using the technique described above to
 the moments obtained by the direct numerical simulation.
Note that for the reconstructed probability distributions of the process we  have two sources of error:
the approximation error of the moment closure and the error associated with the maximum entropy
reconstruction as explained below. 
In our experimental results
we therefore apply the reconstruction to both moments obtained from
the moment equations as well as the more accurate approximation obtained from a direct numerical simulation.
The latter, however, is only possible for systems where the average molecule numbers remain small since otherwise too many states have to be considered during the integration.

\section{Maximum Entropy Reconstruction}\label{sec:maxent}
	The moment closure is usually used to approximate the moments
	of a stochastic dynamical system over time.
	The numerical integration of the correspondent ODE system is
	usually faster than a direct integration of the probability distribution
	or an estimation of the moments based on  Monte-Carlo simulations of the system.
	However, if one is interested in certain events and only the moments of the distribution are known,
	the corresponding probabilities are not directly accessible and
	have to be reconstructed based on the moments.
	Here, we shortly review standard approaches to reconstruct 
	one-dimensional 
	marginal probability distributions
	$\pi_i(x_i,t) = P(X_i(t)=x_i)$ of a Markov chain that describes 
	the dynamics of chemical reactions network.
	The task of approximating multi-dimensional distributions 
	follows the	same line however for our case these techniques 
	revealed to be not effective due to numerical difficulties
	in the optimization procedure.
	Thus, we have given (an approximation of) the moments of the $i$-th population and obviously,
	the corresponding distribution is in general not uniquely determined for a finite set of moments.
	In order to select one distribution from this set, we apply the maximum entropy principle.
	In this way we minimize the amount of prior information about the distribution 
	and avoid any other latent assumption about the distribution.
	Taking its roots in statistical mechanics and thermodynamics~\cite{PhysRev.106.620}, 
	the maximum entropy approach was successfully applied 
	to solve moment problems in the	field of 
	climate prediction~\cite{abramov2005information,kleeman2002measuring,roulston2002evaluating},
	econometrics~\cite{wu_calculation_2003},
	performance analysis~\cite{tari_unified_2005,guiasu1986maximum}
	 the and many others.

	\subsection{Maximum Entropy Approach}
	The maximum entropy principle says that
	among the set of allowed discrete probability distributions $\distrFamily$
 	we choose the probability distribution $q$ that maximizes the entropy $H(g)$ over all  distributions $g \in \distrFamily$, i.e.,
	\begin{equation}
	\label{eq:maxShannonProblem}
	\begin{array}{c}
		\maxent = \arg \max_{g\in \distrFamily} H(g)
				= \arg \max_{g\in \distrFamily} \left( -\sum_x g(x) \ln{g(x)}\right)
				.
	\end{array}
	\end{equation}
	where $x$ ranges over all possible states of the discrete state space.
	Note that we assume that all distributions are defined on the same state space.	In our case the set $\distrFamily$
	consists of all discrete probability distributions that satisfy the moment constraints. Given a sequence of 	$M$ non-central moments
	$$\Ex{X^k}=\mu_k, k=0,1,\ldots,M,$$
	 the following constraints are considered
	\begin{equation}\label{eq:momentconstr}
	\sum_x x^k g(x)   = \mu_k, k=0,1,\ldots,M.
	\end{equation}
	Here, we choose $g$ to be a non-negative function 
	and add the constraint $\mu_0=1$ in order 
	to ensure that $g$ is a distribution.
	The above problem is a nonlinear constrained optimization problem, which is usually
	addressed by the method of Lagrange. Consider the Lagrangian
	functional
	\begin{equation*}
	\begin{array}{c}
		\mathcal{L}(g,\lambda)
		= H(g) - \sum\limits_{k=0}^{M} \lambda_k
				\left( \sum_x x^k g(x)  - \mu_k \right),
	\end{array}
	\end{equation*}
	where $\lambda=(\lambda_0,\ldots,\lambda_M)$ are the corresponding Lagrangian multipliers.
	It is possible to show that maximizing the  unconstraint Lagrangian $\mathcal{L}$
	gives a solution to the constrained maximum entropy problem
 	The variation of the functional $\mathcal{L}$ according to
	the unknown distribution provides the general form of $g(x)$
	$$	\frac{\partial \mathcal{L}}{\partial g(x)} = 0
		\implies
		g(x) = \exp \left( -1 -\sum\limits_{k=0}^{M} \lambda_k x^k \right)
		=\frac{1}{Z(x)} \exp \left( -\sum\limits_{k=1}^{M} \lambda_k x^k \right),
	$$
	where
	$$Z(x) = e^{1+\lambda_0}
	      = \displaystyle\sum_x \exp \left( -\sum\limits_{k=1}^{M} \lambda_k x^k \right) $$
	is a normalization constant.
	In 
	the dual approach
	we insert the above equation for $g(x)$ into
	the Lagrangian thus we can transform the problem into an
	unconstrained convex minimization problem of the dual function w.r.t
	to the dual variable $\lambda$
	$$\Psi(\lambda)=\ln Z(x) + \sum\limits_{k=1}^{M} \lambda_k \mu_k,$$
		According to the Kuhn-Tucker theorem, the solution
	$\lambda^* = \arg \min \Psi(\lambda)$ of the minimization problem
	for the dual function equals the solution $q$ of the original constrained optimization problem~\eqref{eq:maxShannonProblem}.

	\subsection{Maximum Entropy Numerical Approximation}
	It is possible to solve  the constrained maximization problem~\eqref{eq:maxShannonProblem}
	for $m \leq 2$ analytically. For $m>2$ 
	  numerical methods have to be applied
	to incorporate the knowledge of  moments of order three and more.
 	Possible numerical solution techniques include the Newton minimization procedure~\cite{mead1984maximum}, the
	iterative minimization~\cite{bandyopadhyay2005maximum}
	and the application of the Broyden-Fletcher-Goldfarb-Shanno (BFGS) procedure~\cite{byrd1995limited}.
	For our experimental results, we used the  algorithm proposed by 
	Abramov~\cite{abramov2010multidimensional} and the corresponding software~\footnote[1]{The software is 
	available at
	http://homepages.math.uic.edu/$\sim$abramov/}.
	It can be used to reconstruct distributions 
	with up to 4 dimensions.
	Due to the difficulties in the numerical optimization
	procedure here we restrict ourselves to the 
	the reconstruction of one-dimensional marginal distributions 
	$q_i(\cdot)$ that approximate $\pi_i(\cdot,t)$.
	The essential idea of the algorithm is 
	preconditioning of the
	original problem to overcome numerical difficulties.
	This standardization is conducted through a sequence of
	linear transformations of coordinates in  phase space.
	They include shifting of the moments,
	rescaling its values to shorten the difference between orders of magnitudes
	and representing the optimization problem in the generalized orthogonal polynomial basis
	to lower the sensitivity of the Lagrange multipliers to high values
	of monomials $x^i$.
	After this sequence of preconditioning steps is completed,
	the optimization BFGS procedure is used to obtain the
	Lagrangian multipliers $\lambda^*$.
	Solving the ODE system for moments we obtain an approximation 
	of the moment values and thus the distribution with such moments
	might not exist.
	This is not a problem in our case since 
	both the number of points in the Gauss-Hermite
	quadrature formula used in integral computation
	and the tolerance for optimization procedure are 
	controllable in the implementation we use,
	so that the optima (if it exists) can be found
	even using the approximated values.
	Theoretical conditions for existence of the solution for moment problem
	are elaborated in detail in~\cite{tari_unified_2005,Stoyanov_2000,Lin199785}.
	The similar analysis for the multivariate case is provided in~\cite{kleiber_multivariate_2013}.

	Abramov focuses on continuous state spaces and reconstructs densities instead of discrete
	probabilities. However, it turns out that the computed densities can be easily transformed into
	discrete distributions (cf. Section~\ref{sec:exp}).
	The probability distributions $\pi_i(\cdot,t)$ we are interested in are discrete
	and in order to approximate the discrete probability $P \left( X_i = x_i \right)$
	for a molecule count $x_i$ of the $i$-th species
	we integrate the reconstructed density $q_i$ and set
	\begin{equation}
	\label{eq:aux_discrete_rv}
	\begin{array}{c}
		\tilde\pi_i(x_i,t) = 
		\begin{cases}
		 \int\limits_{x_i - \frac{1}{2}}^{x_i + \frac{1}{2}} q_i(x) \, \mathrm{dx}, & x_i > 0 \\
		2  \cdot \int\limits_{0}^{x_i + \frac{1}{2}} q_i(x) \, \mathrm{dx}, & x_i=0
		\end{cases}
	\end{array}
	\end{equation}
		Additionally we have to  truncate the support ${\mathbb R}$ of $q_i$ to ${\mathbb R^+}$ such that
	$q_i(x) = 0$ for all $x \in (-\infty,0)$. This can be done by 
	adjusting the normalization constant $Z(x)$ accordingly.
	The two of three models we are considering have the unbounded state space
	which does not allow to simplify the process of moment computation in the
	maximum entropy optimization procedure thus they are approximated 
	by Gauss-Hermite quadrature formula.

\begin{table}[b]
\tbl{Errors of the moment closure approximation for the 
dimerization
network\label{tab:dimerization_mc_error}}{
\begin{tabular}{|c|c|c|c|c|c|c|}
\hline
ord. & $\#$ equ. & error ord. 1& error ord. 2 & error ord. 3 & error ord. 4 & error ord. 5\\
\hline
2 & 5 & 0.001754 & 0.003495 &- &-  & -\\
\hline
3 & 9 & 0.001752 & 0.003492 & 0.005215 &- & -\\
\hline
4 & 14 & 0.001743 & 0.003465 & 0.005211 & 0.006907 & -\\
\hline
5 &  20 & 0.001721 & 0.003418 & 0.005183 & 0.006901 & 0.008555 \\
\hline
\end{tabular}}
\end{table}

\section{Numerical Results}\label{sec:exp}
In this section we show experimental results of the maximum entropy approach when it is applied to the 
moments of a reaction network. To estimate the quality of the probability distribution reconstruction,
we compare the obtained distributions to those obtained  via a direct numerical simulation. Thus, we 
consider only systems which are small and where a direct numerical simulation is possible. 
Clearly, for more complex systems with high population sizes
a direct numerical simulation is not feasible while the  running time of the 
moment closure approximation is independent of the population sizes.

In order to distinguish errors that are introduced by the moment closure approximation
from errors     introduced by the reconstruction, we also compare the obtained moments
with those computed based on the distributions obtained via direct numerical simulation. 
Moreover, we applied the maximum entropy approach to the 
more accurate approximation of the moments obtained via direct numerical simulation.
It is important to note that 
in all cases the maximum entropy optimization only takes less than one second.
We therefore only list the  running  time  of the moment closure method.
 \subsection{Simple Dimerization system}
We first consider Example~\ref{ex:simple_dimerization_reactions}
and investigate the numerical accuracy of the moment equations
for this example. The moment closure approximation only takes less than one second 
for this example and
Table~\ref{tab:dimerization_mc_error} compares the moments approximated 
with the moment closure  with the moments obtained from the direct numerical simulation 
(as described in Section~\ref{sec:directnumsim}).

\begin{figure}[t]
\centering
  \subfloat[][]{\includegraphics[width=0.5\textwidth]{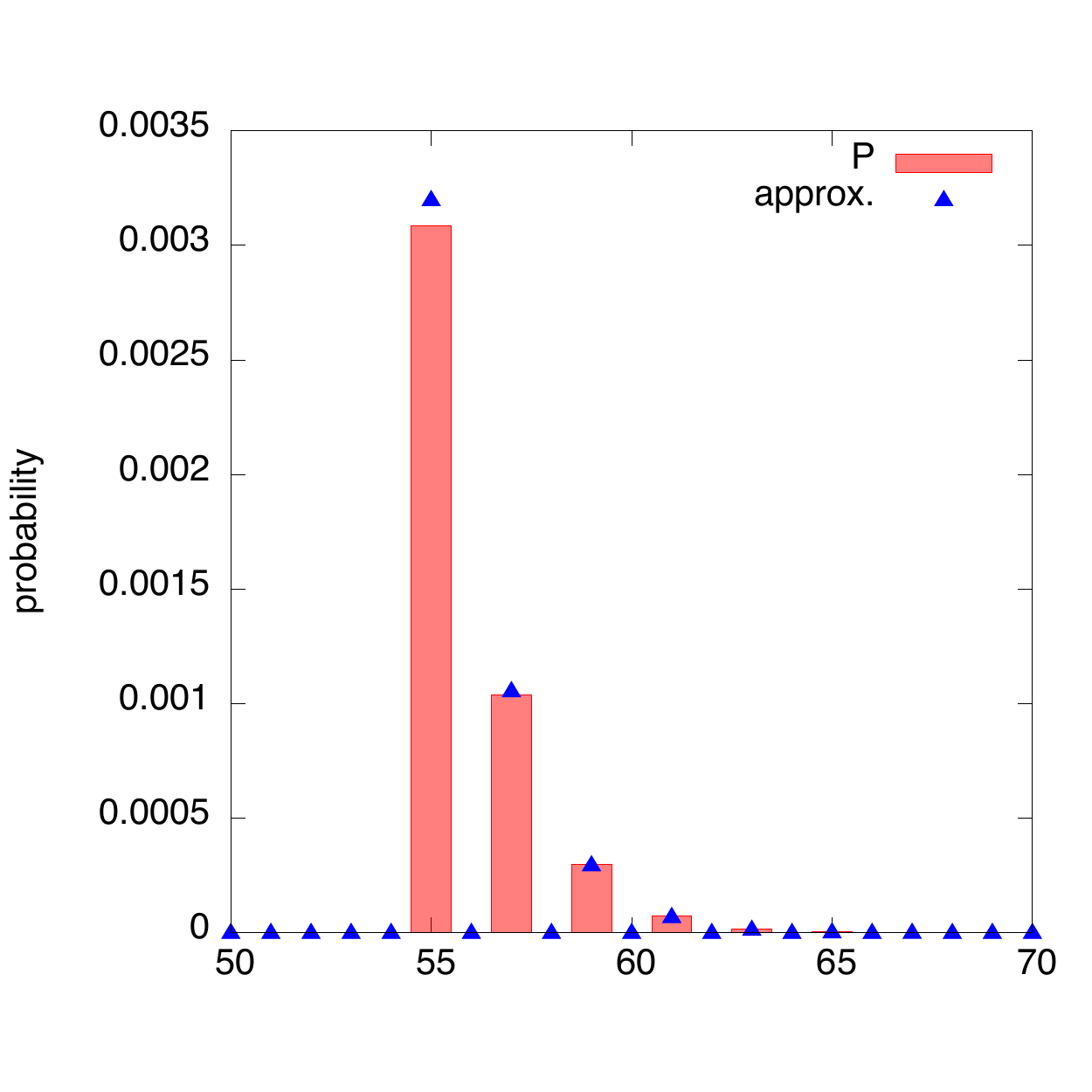}}
  \subfloat[][]{\includegraphics[width=0.5\textwidth]{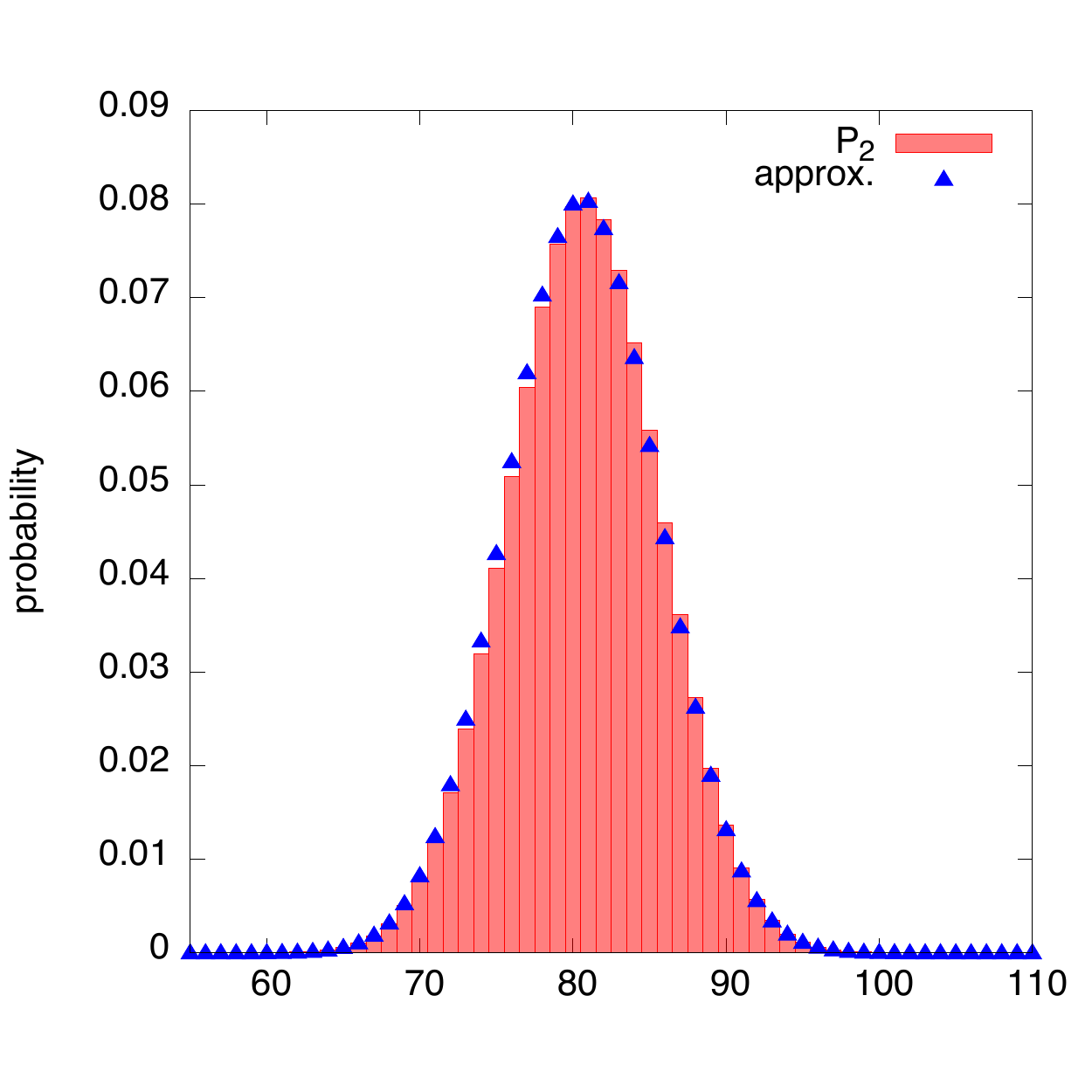}}
	\caption{\label{fig:mer_simple_dimerization_pdfs}
	Maximum entropy reconstruction of marginal probability distributions of the protein counts
    (a) $P$ and (b) $P_2$ at time instant $t=20$ for simple dimerization system.}
\end{figure}

\begin{table}[b]
\tbl{Maximum entropy reconstruction results for the dimerization network\label{tab:dimerization_mer_error}}{
\begin{tabular}{|c|c|c|c|c|}
\hline
ord. & $||\epsilon_P||_\infty$ 
& $||\epsilon^*_P||_\infty$ 
& $||\epsilon_{P_2}||_\infty$ 
& $||\epsilon^*_{P_2}||_\infty$ 
\\
\hline
2 & 0.001764 & 0.000623 & 0.001764 & 0.000136\\
\hline
3 & 0.000860 & 0.000054 & 0.001782 & 0.000623\\
\hline
4 & 0.001683 & 0.000136 & 0.001683 & 0.000053\\
\hline
5 & 0.001641 & 0.000132 & 0.001691 & 0.000136\\
\hline
\end{tabular}}
\end{table}

	The first column 
	refers to the highest order of the moments 
	that was considered during the computation, 
	i.e.,  all central moments of higher order are approximated with zero.  
	In addition we list the relative errors of the means (error ord. 1) and 
	the moments of higher order $k$ using the error norm
	$$\max\limits_{i} \frac{|\hat{m}_i^{(k)}-m_i^{(k)}|}{m_i^{(k)}}.$$
	Here $\hat{m}_i^{(k)}$ and $m_i^{(k)}$
	are the values of moments $E(X_i^k)$
	computed using the moment closure method and obtained via direct numerical simulation
	and the maximum is taken over the 
	chemical species.
	The second column in tables refers to the number of equations that were
	integrated for the moment closure method.
	The initial protein numbers are chosen as $P=301$ and $P_2=0$ and we consider the
	system  at time $t=20$~\cite{StumpfJournal}.
	We find that the moment closure approximation provides very accurate results.


 Next, we reconstruct the marginal distributions of the species $P$ and $P_2$  and compare then with those
obtained using the direct numerical simulation  (where we chose $\delta=10^{-15}$ yielding a total approximation error of $\epsilon=5\cdot 10^{-15}$, see also Eq.~\eqref{eq:eps}).
For instance, to reconstruct the distribution of   $P$
we used the sequence of moments $\mu_0,\ldots,\mu_k$, where
$\mu_j = E(X^j_1(t))$ and  $X_1(t)$ represents the number of molecules of type $P$ at time $t$. The values of moments $\mu_0,\ldots,\mu_k$ 
are approximated by the moment closure method. 
In Fig.~\ref{fig:mer_simple_dimerization_pdfs} we plot the 
distribution of $P$ and $P_2$ where we use red bars for the 
distribution obtained via direct numerical simulation and  blue
triangles  for the reconstructed distribution (moment closure approximation and maximum entropy reconstruction).
For the order of the moments that were considered, we used in 
for both species the order with which the best approximation was
obtained (see below). 
 
In Table~\ref{tab:dimerization_mer_error} we show how accurate 
  the approximation of individual probabilities $\pi_i(x,t)$ is
by calculating the Chebyshev 
distance
$$||\epsilon_{i}||_\infty = \max\limits_{x} | \pi_i(x,t) - \tilde\pi_i(x,t)|,$$
where $\pi_i(x,t)$ is the "true" probability 
of having $x$ molecules of type $i$ at time $t$
(obtained via the accurate direct numerical simulation) 
and $\tilde \pi_i(x,t)$ is the value obtained from 
the combination of moment closure approximation and maximum entropy reconstruction.
The distance is calculated for all non-negative integers $x$ for which
  $\pi_i(x,t)$   is non-negative,
i.e. for species $P$ we only consider every second integer value.
In particular, we approximate the marginal distribution for $P$ 
by using a modified form of Eq.~\ref{eq:aux_discrete_rv}:
	\begin{equation*}
	\begin{array}{c}
		\tilde\pi_i(x_i,t) = 
		\begin{cases}
		 \int\limits_{x_i-1}^{x_i+1} q_i(x) \, \mathrm{dx}, & x_i > 0 \\
		2  \cdot \int\limits_{0}^{x_i+1} q_i(x) \, \mathrm{dx}, & x_i=0
		\end{cases}
	\end{array}
	\end{equation*}
\begin{table}[b]
\tbl{Moment closure approximation results for the multi-attractor network\label{tab:multiattractor_mc_error}}{
\begin{tabular}{|c|c|c|c|c|c|c|c|}
\hline
ord. & $\#$ equ. &  time (sec) & error ord. 1& error ord. 2 & error ord. 3 & error ord. 4 & error ord. 5\\
\hline
2 & 104 & 2 & 0.043987 & 0.077507 &- &-  & -\\
\hline
3 & 559 & 40 & 0.043987 & 0.067790 & 0.104288 &- & -\\
\hline
4 & 2379 & 443 & 0.043987 & 0.058938 & 0.082293 & 0.096345 & -\\
\hline
5 & 8567 & 3649 & 0.043987 & 0.037542 & 0.066227 & 0.056258 & 0.110358 \\
\hline
\end{tabular}}
\end{table}
In addition, we give the error~$||\epsilon^*_i||_\infty$ for the case where the maximum entropy 
reconstruction was applied to the moments calculated from the results of the direct numerical simulation.
We observe that the best maximum entropy method provides
the least error when the moments of order up to 3 (for $P$) and up to 4 (for $P_2$) are used to reconstruct 
the marginal distribution.
However, the reconstruction is very accurate in all cases
and the reason why the Chebyshev distance does not decrease when 
more moments are considered might be that the properties of the
distribution are already well captured by the moments of an order up to three.

\subsection{Multi-attractor system}
\newcommand{\paxp}{\ensuremath{\mathit{PaxProt}}}
\newcommand{\mafap}{\ensuremath{\mathit{MAFAProt}}}
\newcommand{\deltap}{\ensuremath{\mathit{DeltaProt}}}
\newcommand{\paxd}{\ensuremath{\mathit{PaxDna}}}
\newcommand{\mafad}{\ensuremath{\mathit{MAFADna}}}
\newcommand{\deltad}{\ensuremath{\mathit{DeltaDna}}}
Our second case study is the so-called multi-attractor model~\cite{zhou_predicting_2011}.
It consists of 23 chemical reactions (listed in Appendix~\ref{app:A}) and 
describes the dynamics of three genes and the corresponding proteins. 
The proteins 
$\paxp$, $\mafap$ and $\deltap$
are able to bind to the promotor regions of the DNA and
activate or suppress the production of other proteins.
The model is infinite in three dimensions.

Again, we first consider the accuracy of the moment closure approximation (cf. Table~\ref{tab:multiattractor_mc_error}) in the same way as for the previous example but list
the running time in addition (third column). 
The values or stochastic reaction constants are chosen as
$c_p = 5, c_d = 0.1, c_b = 1.0, c_u = 1.0$ and we consider the system at time $t=10$.
As   initial conditions we assumed one molecule for all
DNA-like species   ($\#\paxd = 1$, $\#\mafad = 1$, $\#\deltad = 1$)
and the molecular counts for the remaining species are 0.

We find that
the   moments obtained via the moment closure approximation are accurately approximated.
For instance
the average number of $\mafad$ is approximated as $19.719$
while the result of the direct numerical simulation gives $19.544$.
Note that it takes $20634$ seconds to finish the numerical simulation
(the size of the truncated state space $|S|=7736339$)
whereas the moment closure approximation takes only $3649$ seconds.


%

\begin{table}[b]
\tbl{Maximum entropy reconstruction results for the multi-attractor system
\label{tab:multi_attractor_mer_error_t10}}{
\begin{tabular}{|c|c|c|c|c|c|c|}
\hline
ord. 
& $||\epsilon_{\mbox{\tiny MAFAProt}}||_\infty$ & $||\epsilon^*_{\mbox{\tiny MAFAProt}}||_\infty$ 
& $||\epsilon_{\mbox{\tiny DeltaProt}}||_\infty$ & $||\epsilon^*_{\mbox{\tiny DeltaProt}}||_\infty$
& $||\epsilon_{\mbox{\tiny PaxProt}}||_\infty$ & $||\epsilon^*_{\mbox{\tiny PaxProt}}||_\infty$ 
\\
\hline
2 
&0.016082 & 0.012322 &0.009978 &0.009186 & 0.053345 & 0.033673\\
\hline
3 &0.012075 &0.009372 & 0.009904 &0.009728	& 0.037701 & 0.031701\\
\hline
4 &0.009030 &0.008892 & 0.009541&0.007568	& 0.033590 &0.027062 \\
\hline
5 &0.006117 &0.005308 & 0.007053 &0.005295	& 0.030783 &0.021304  \\
\hline
\end{tabular}}
\end{table}

\begin{figure}[t]
\centering
  \subfloat[][]{\includegraphics[width=0.33\textwidth]{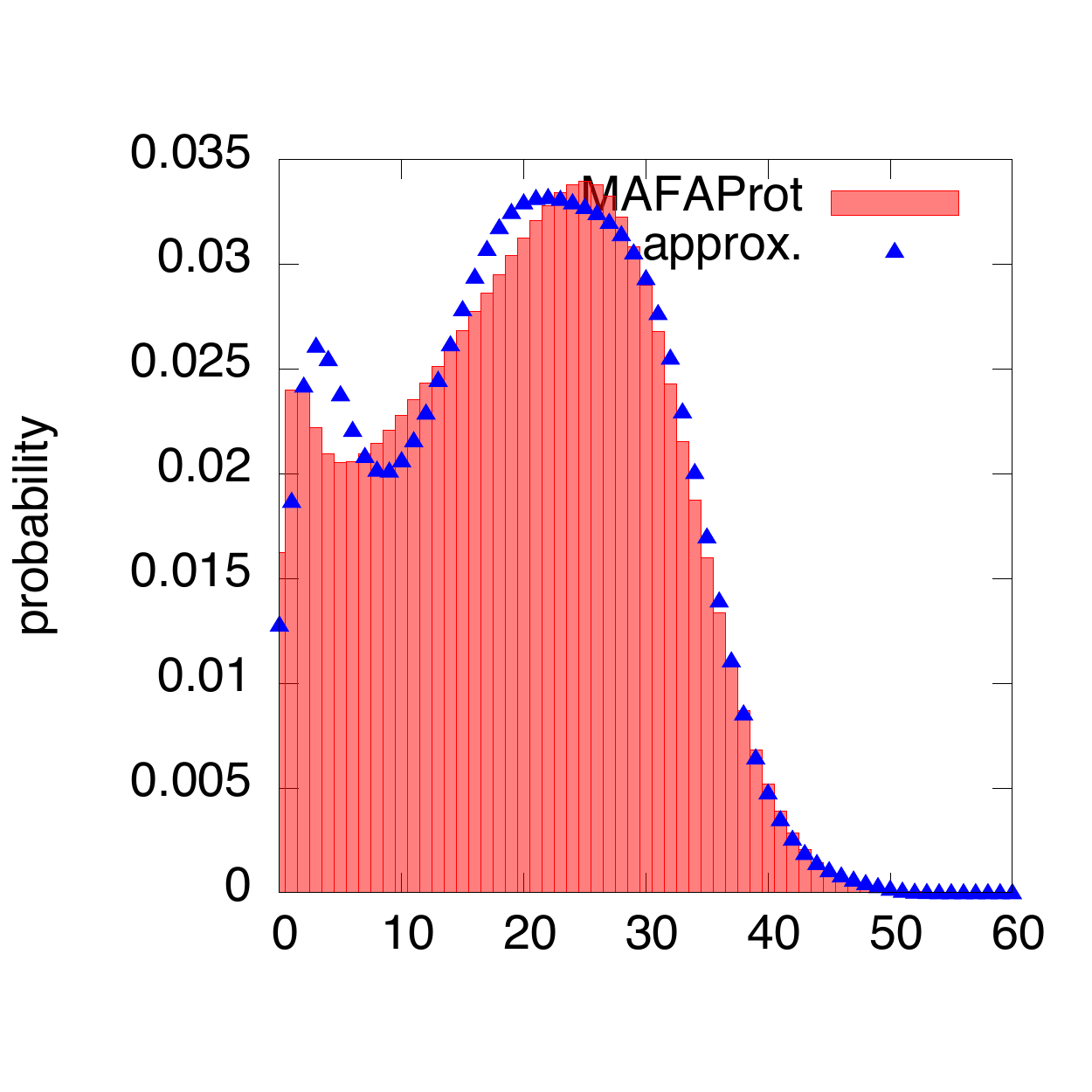}}
  \subfloat[][]{\includegraphics[width=0.33\textwidth]{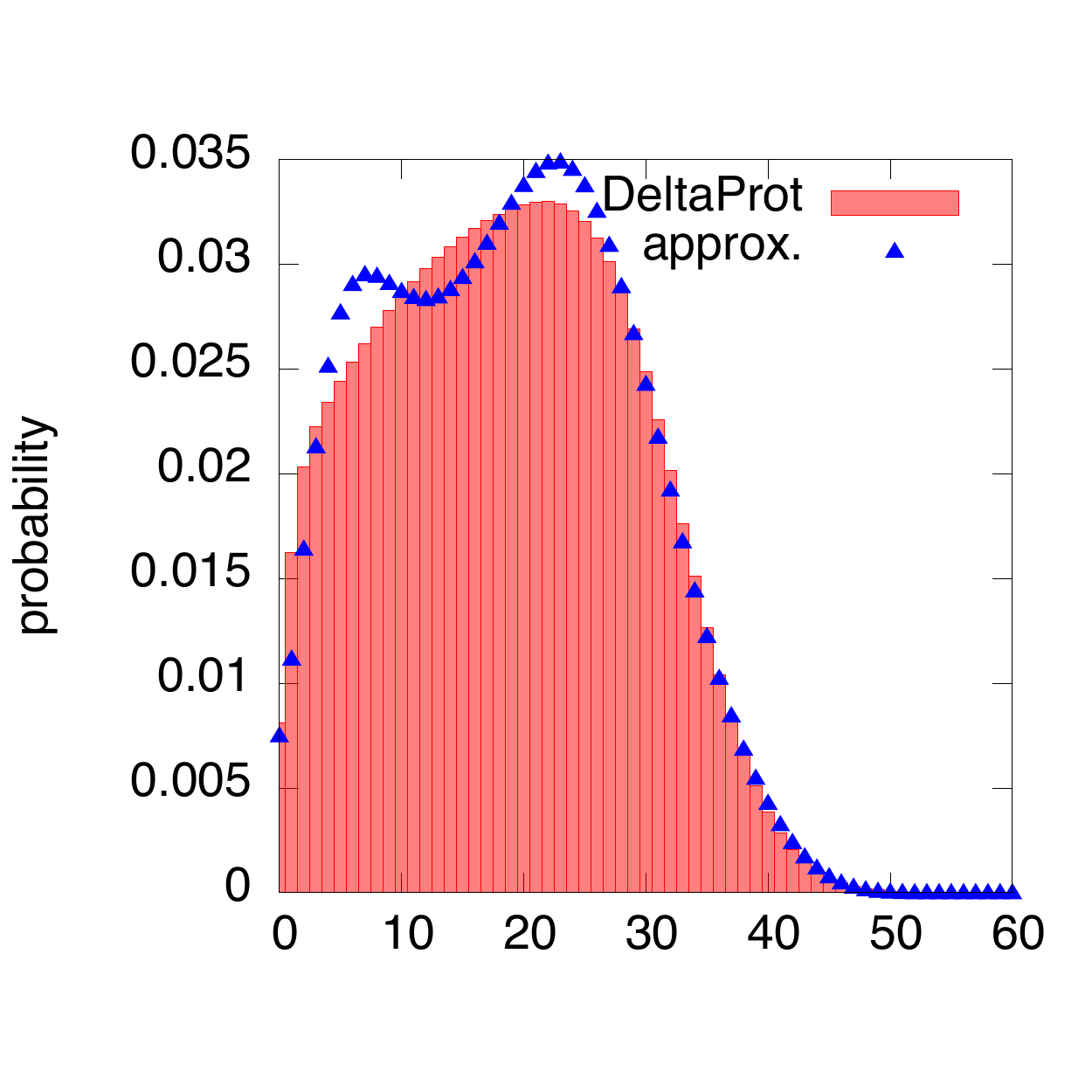}}
  \subfloat[][]{\includegraphics[width=0.33\textwidth]{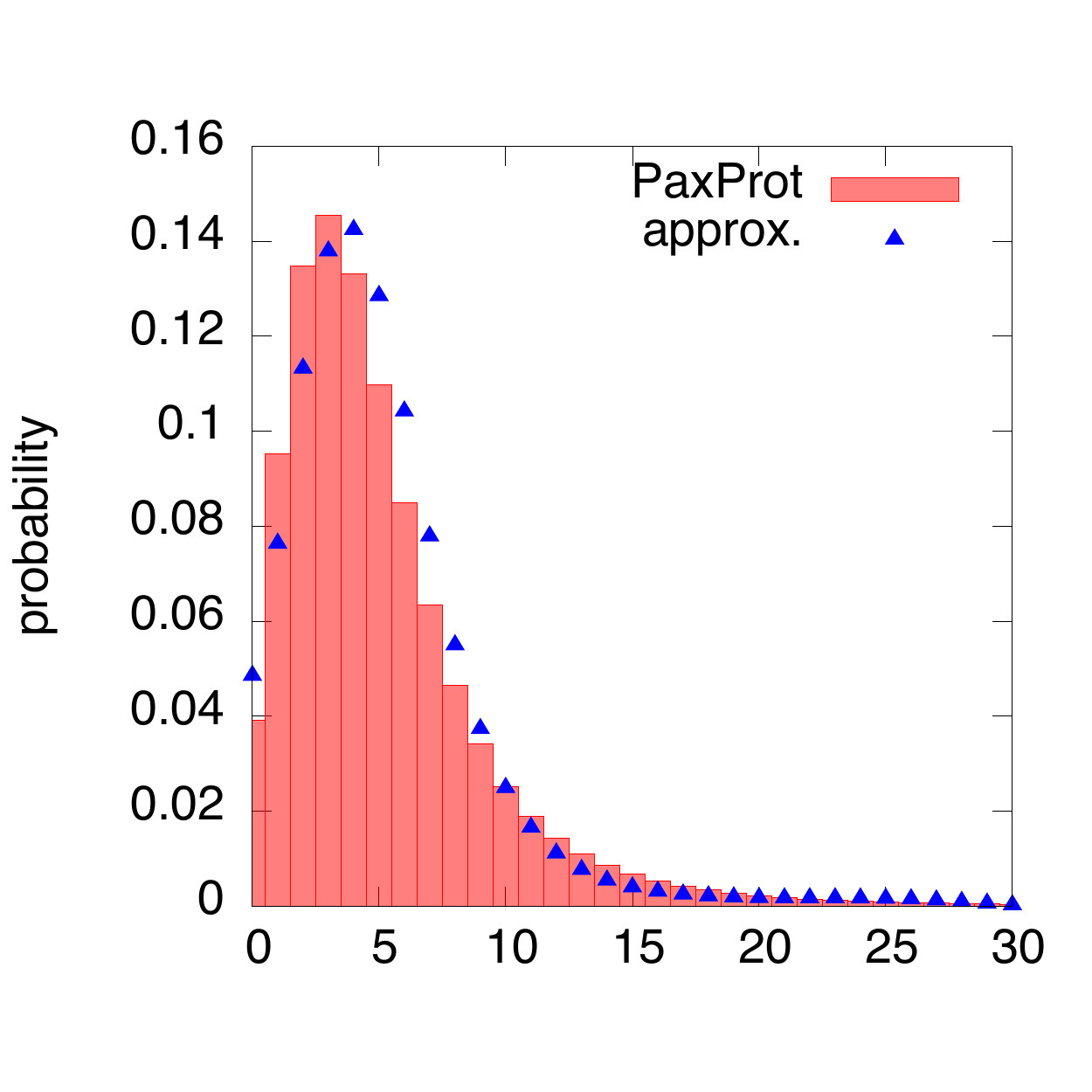}}
	\caption{\label{fig:multi_attractor_t10_pdfs}
	Maximum entropy reconstruction of marginal probability distributions 
	for multi-attractor system at time instant $t=10$.}
\end{figure}

Next we consider the reconstruction of  the marginal distribution of 
$\paxp$, $\mafap$ and $\deltap$. 
The results  are given in  
Table~\ref{tab:multi_attractor_mer_error_t10} and
the best obtained reconstructions are plotted 
in Figure~\ref{fig:multi_attractor_t10_pdfs}
for all three proteins. We compare the results with the solution of the direct numerical simulation
(where we chose $\delta=10^{-15}$ yielding a total approximation error of $\epsilon=6\cdot 10^{-10}$, 
see also Eq.~\eqref{eq:eps}).
We see that the error decreases when
more moments are considered. In particular, if all moments up to order 5 are considered, the error is about an order of magnitude
lower than the state probabilities around the average molecule
count.
 However, the reconstructed distribution of $\deltap$ shows some 
 artifacts,
which are due to the bimodal form of the exponential function
where the power is given by the polynomial of degree five.
We do not observe the similar type of artifacts 
in the reconstructed distribution of $\paxp$ since there
the reconstructed coefficients corresponding to high powers 
are close to $0$.

\subsection{Exclusive switch}
Finally, we consider the exclusive switch
system.
We chose reaction rate constants and initial conditions
  as   in Example~\ref{ex:exSwitchReactions}.
 Again, we first consider the accuracy of the moment closure approximation (cf. Table~\ref{tab:exswitch_mc_error}) at time $t=100$ in the same way as for the previous examples.
 As also noted by Grima, the error of the moments, that have the highest order considered during the computation,  is  rather high~\cite{grima2012study}. Thus, in the moment closure approximation we have to consider
 at least all moment up to order five to accurately estimate the moments up to order four.

\begin{table}[t]
\tbl{Moment closure approximation results for the exclusive switch\label{tab:exswitch_mc_error}}{
\begin{tabular}{|c|c|c|c|c|c|c|c|}
\hline
ord. & $\#$ equ. &  time (sec) & error ord. 1& error ord. 2 & error ord. 3 & error ord. 4 & error ord. 5\\
\hline
2 & 20 & $<1$ & 0.004555 & 0.194240 &- &-  & -\\
\hline
3 & 55 & $<1$ & 0.004555 & 0.026281 & 0.060490 &- & -\\
\hline
4 & 125 & 2 & 0.004555 & 0.020493 & 0.028242 & 0.136965 & -\\
\hline
5 &  251 & 6 & 0.004555 & 0.017774 & 0.027933 & 0.026724 & 0.015461 \\
\hline
\end{tabular}}
\end{table}

\begin{figure}[t]
\centering
  \subfloat[][]{\includegraphics[width=0.5\textwidth]{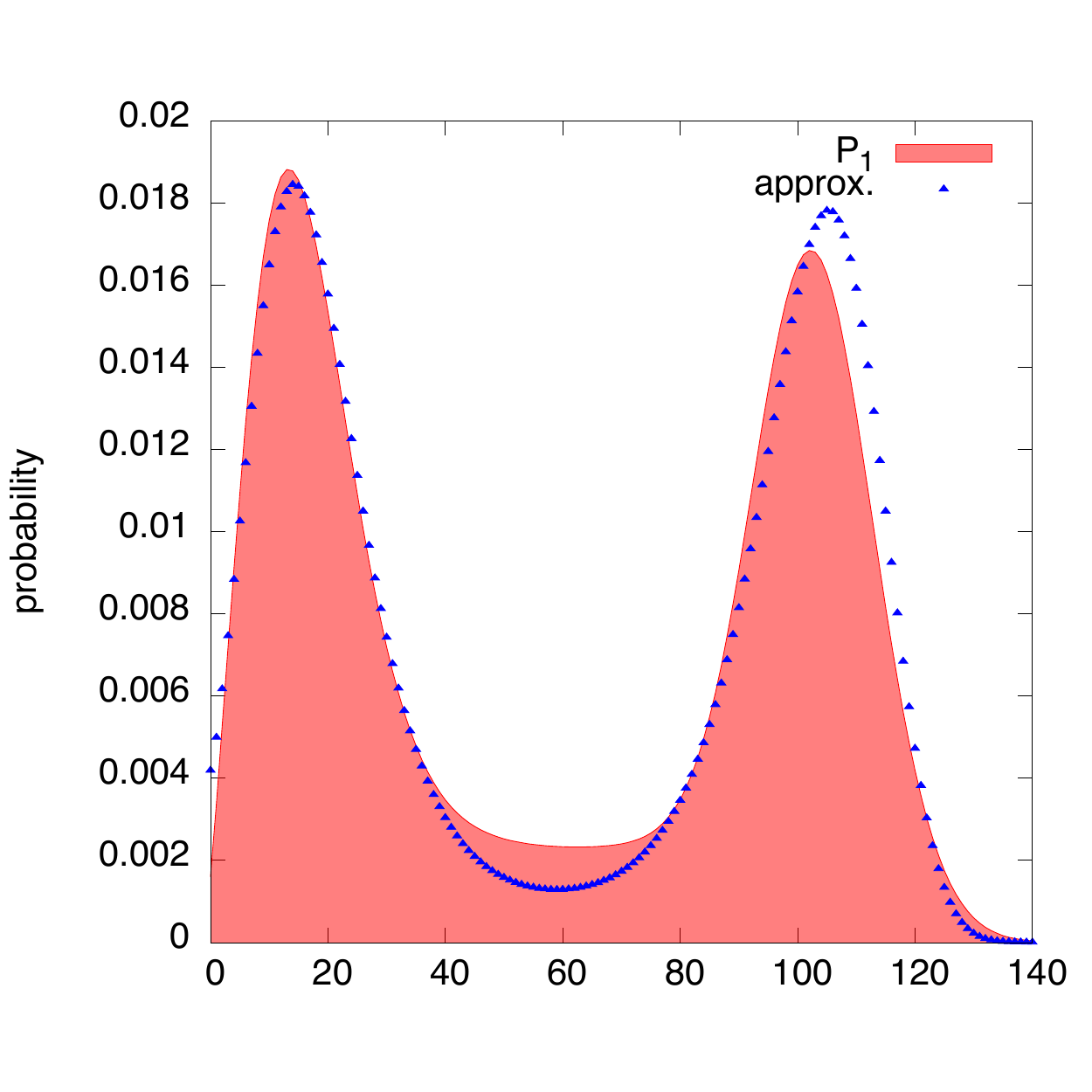}}
  \subfloat[][]{\includegraphics[width=0.5\textwidth]{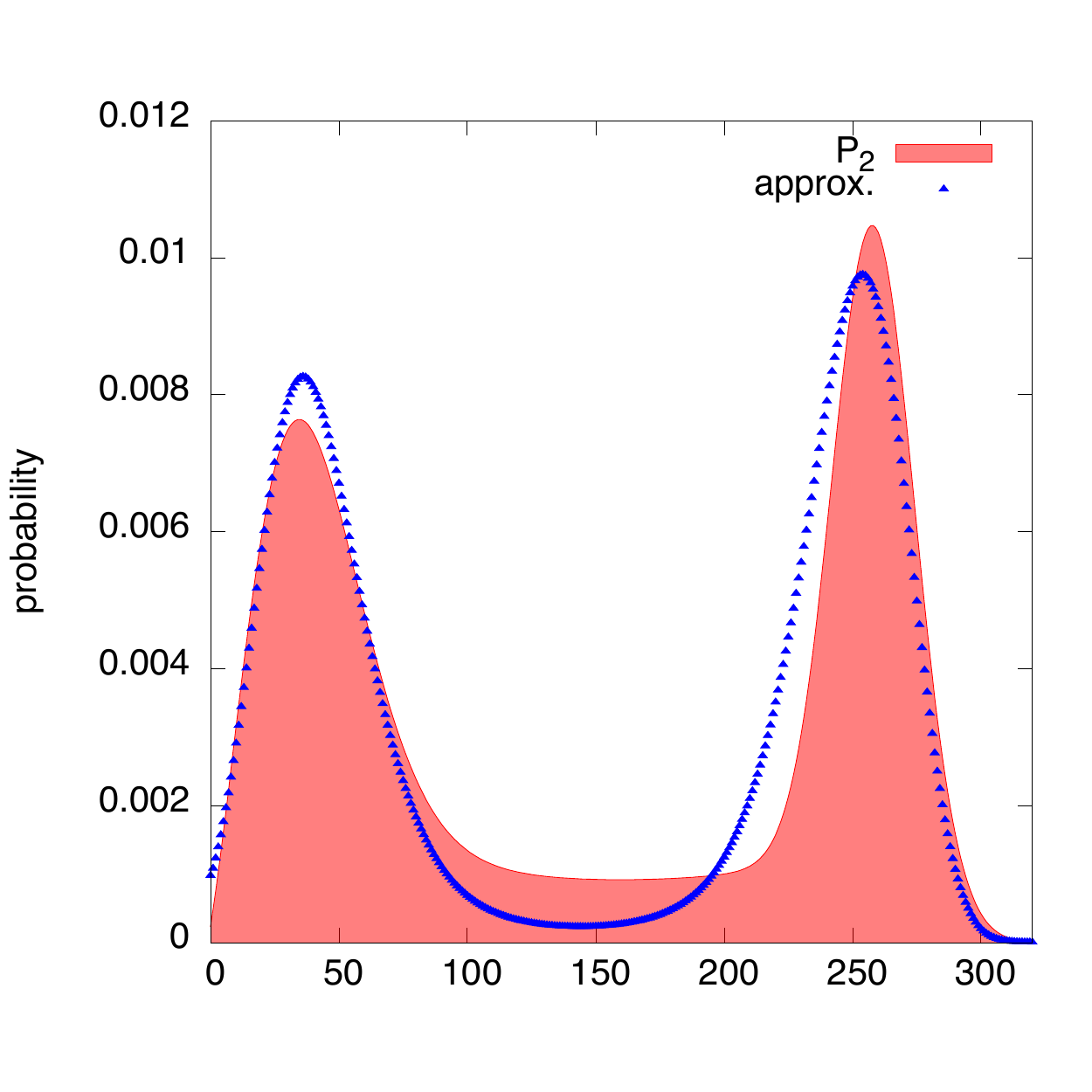}}
	\caption{\label{fig:mer_exSwitch_t60_pdfs}
	Maximum entropy reconstruction of marginal probability distributions of the protein counts
    (a) $P_1$ and (b) $P_2$ at time instant $t=60$ for exclusive switch system.}
\end{figure}

\begin{figure}[t]
\centering
  \subfloat[][]{\includegraphics[width=0.5\textwidth]{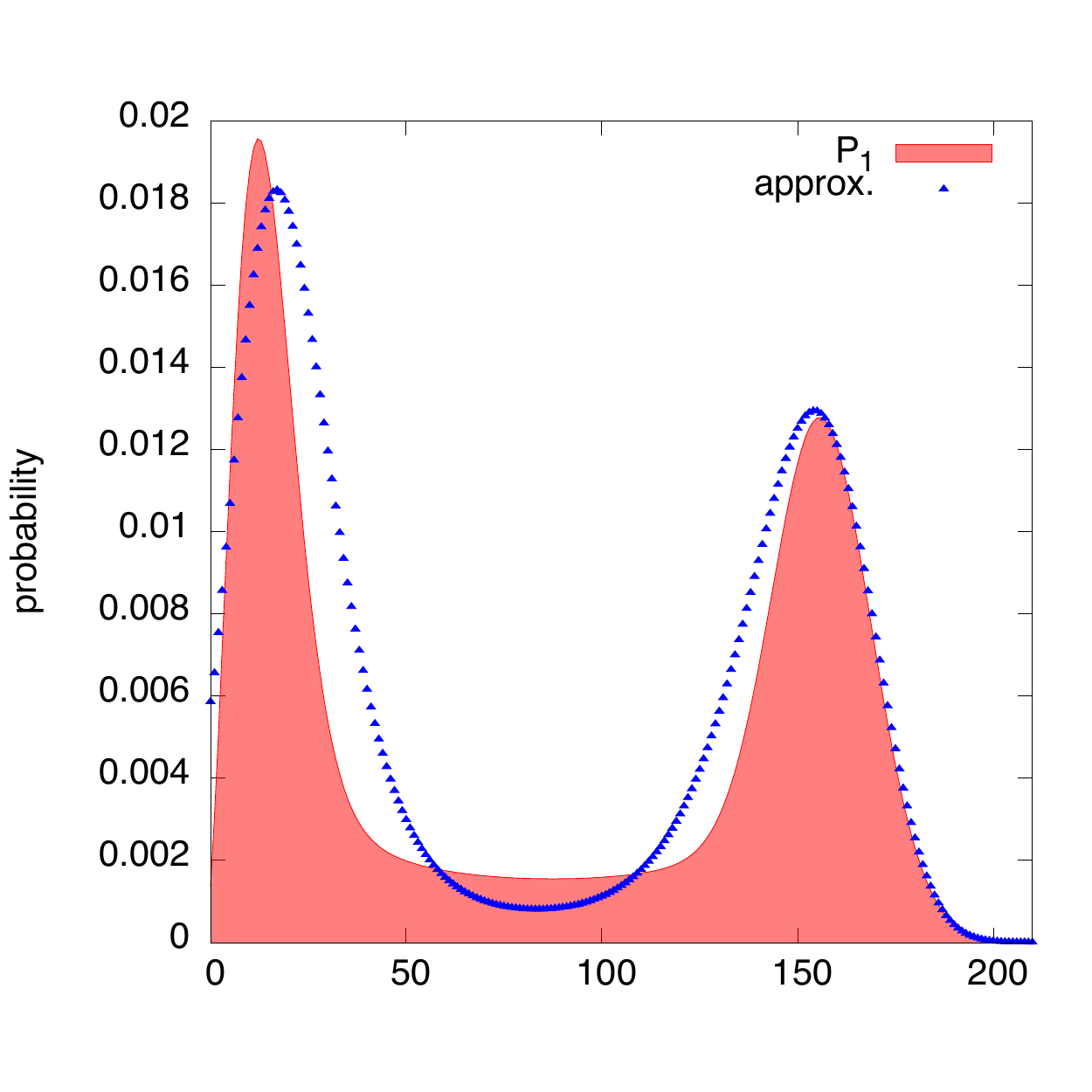}}
  \subfloat[][]{\includegraphics[width=0.5\textwidth]{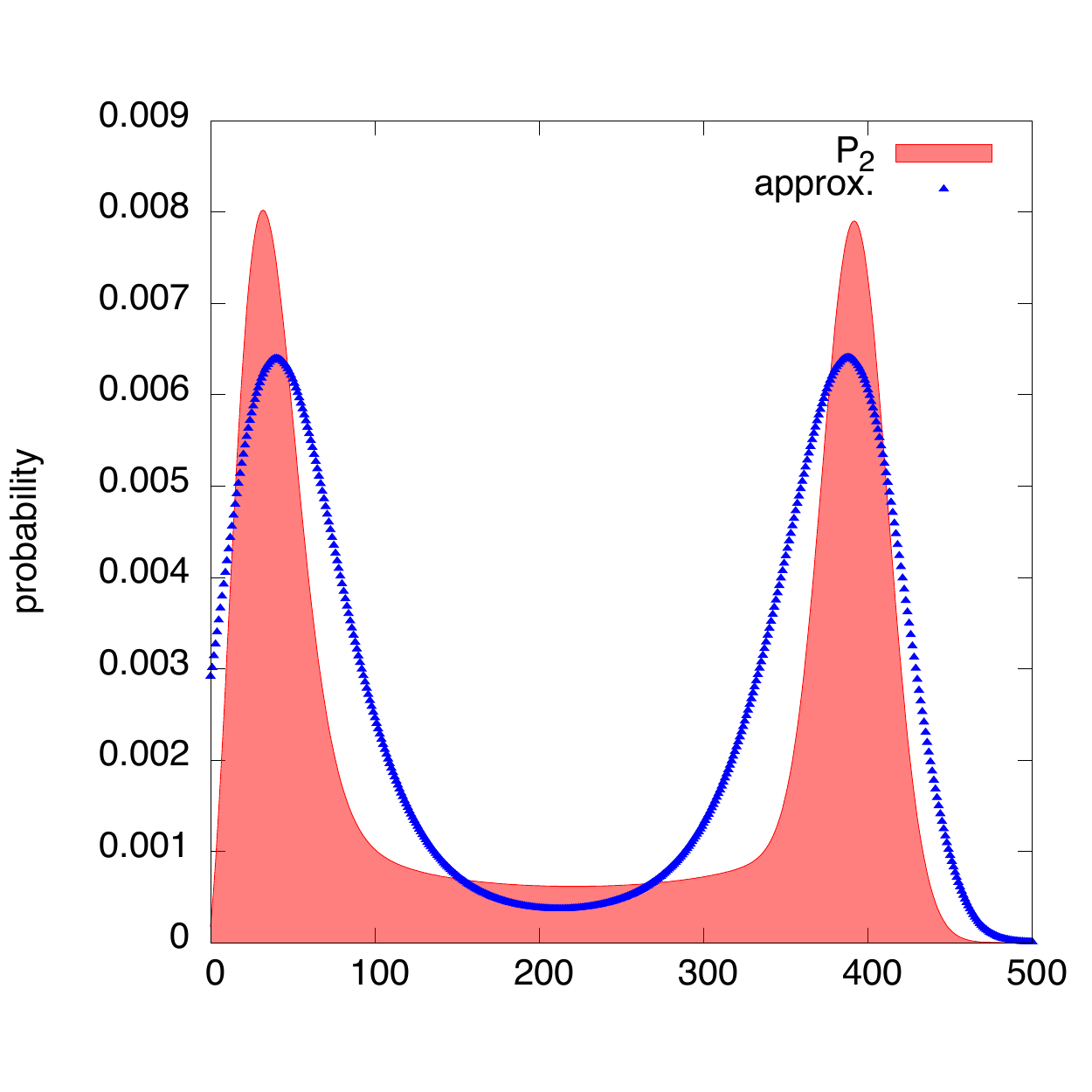}}
	\caption{\label{fig:mer_exSwitch_t100_pdfs}
	Maximum entropy reconstruction of marginal probability distributions of the protein counts
    (a) $P_1$ and (b) $P_2$ at time instant $t=100$ for exclusive switch system.}
\end{figure}

Next we compare  the  marginal distributions of proteins
$P_1$ and $P_2$ at $t=60$ (cf. Table~\ref{tab:exswitch_mer_error_t60})
and at $t=100$ (cf. Table~\ref{tab:exswitch_mer_error_t100})
obtained via a direct numerical simulation (choosing $\delta=10^{-15}$ yielding $\epsilon=6\cdot 10^{-11}$ at time $t=60$ and $\epsilon=8\cdot 10^{-11}$ at time $t=100$)
with the distributions obtained from the maximum entropy reconstruction.
We see that the qualitative property of the system, the bimodality, is well-described
by the   moments up to an order of at least four. 
Thus, it is possible to encode such qualitative properties in the moments. 
The corresponding plots of the marginal distributions and their 
reconstructions are given in Figure~\ref{fig:mer_exSwitch_t60_pdfs} for $t=60$
(Figure~\ref{fig:mer_exSwitch_t100_pdfs} for $t=100$).
If only the means and covariances
are considered, the distribution is not accurately reconstructed,
 $||\epsilon^*_{P_1}||_\infty$ is of the same order as the maximal state probabilities.
As expected the error decreases when moments of higher order are taken into account.




\begin{table}[b]
\tbl{Maximum entropy reconstruction results for the exclusive switch\label{tab:exswitch_mer_error_t60} at time $t=60$}{
\begin{tabular}{|c|c|c|c|c|}
\hline
ord. 
& $||\epsilon_{P_1}||_\infty$ & $||\epsilon^*_{P_1}||_\infty$ 
& $||\epsilon_{P_2}||_\infty$ & $||\epsilon^*_{P_2}||_\infty$ 
\\
\hline
2 & 0.013655 & 0.013630 & 0.008134 & 0.0081149 \\
\hline
3 & 0.013409 & 0.013206 & 0.008074 & 0.007808 \\
\hline
4 & 0.004844 & 0.003084 & 0.002154 & 0.002117\\
\hline
5 & 0.003837 & 0.002541 & 0.001878 & 0.001757\\
\hline
\end{tabular}}
\end{table}

\begin{table}[b]
\tbl{Maximum entropy reconstruction results for the exclusive switch\label{tab:exswitch_mer_error_t100} at time $t=100$}{
\begin{tabular}{|c|c|c|c|c|}
\hline
ord. 
& $||\epsilon_{P_1}||_\infty$ & $||\epsilon^*_{P_1}||_\infty$
& $||\epsilon_{P_2}||_\infty$ & $||\epsilon^*_{P_2}||_\infty$ 
\\
\hline
2 & 0.016287 & 0.016281 & 0.006732 & 0.006563\\
\hline
3 & 0.016270 & 0.016253 & 0.006746 & 0.005158\\
\hline
4 & 0.007783 & 0.007277 & 0.003983 & 0.003301\\
\hline
5 & 0.007527 & 0.007016 & 0.002733 & 0.002455\\
\hline
\end{tabular}}
\end{table}

\section{Conclusions}\label{sec:conc}
We investigated the accuracy and efficiency of a combination of two methods, the moment closure method and the maximum entropy method, which can be used to 
analyze
stochastic models of biochemical reaction networks. With the chemical master equation as a starting point we described how the moments of the
corresponding probability distributions can   be integrated
efficiently over time and how a distribution can be reconstructed
based on the moments. 
Our experimental results show that the proposed combination
of methods has many advantages. It is a fast and surprisingly accurate way of obtaining the 
distribution of the system at specific points in time and therefore well suited for  computationally
expensive tasks such as the approximation of likelihoods or event probabilities.

As future work, we plan to extend the reconstruction procedure 
in several ways. First, we want to consider moments of higher order than five.
Since in this case the concrete values become very large it might be advantageous 
to consider central moments instead which implies that the reconstruction procedure
has to be adapted.
Alternatively, we might (instead of algebraic moments) consider 
other functions of the random variables such as
exponential functions~\cite{mnatsakanov_note_2013},
Fup functions~\cite{gotovac_maximum_2009},
and
Chebyshev polynomials~\cite{bandyopadhyay2005maximum}.
 Another possible extension could
address  the problem of truncating the support of the distribution
such that the reconstruction is applied to a finite support. 
We expect that in this case the reconstruction will become more
accurate since we will not have to rely on the Gauss-Hermite quadrature formula.
%
For instance, the theory of Christoffel functions~\cite{gavriliadis_truncated_2012}
could be used to determine the region where
the main part of the probability mass is located. 

Finally, we want to improve the approximation for species that are present
in very small quantities, since for those species a direct representation of the probabilities 
is more appropriate than a moment representation. Therefore  we plan to consider 
the conditional moments approach~\cite{MCM_Hasenauer_Wolf}, 
where we only integrate the moments of species having large molecular counts but 
keep the discrete probabilities   for the species with small populations.

%
%


\begin{acks}
This research has been partially funded by  the German Research
     Council (DFG) as part of the Cluster of Excellence
     on Multimodal Computing and Interaction at Saarland University and
     the Transregional Collaborative
     Research Center ``Automatic Verification and Analysis of Complex
     Systems'' (SFB/TR~14~AVACS).
     \end{acks}

\bibliographystyle{ACM-Reference-Format-Journals}
\bibliography{bibfile}

\received{February 2014}{XX 2014}{XX 2014}
\clearpage 

\begin{appendix}

\newcommand{\paxddeltap}{\ensuremath{\mathit{PaxDnaDeltaProt}}}
\newcommand{\mafadpaxp}{\ensuremath{\mathit{MAFADnaPaxProt}}}
\newcommand{\mafadmafap}{\ensuremath{\mathit{MAFADnaMAFAProt}}}
\newcommand{\mafaddeltap}{\ensuremath{\mathit{MAFADnaDeltaProt}}}
\newcommand{\deltadpaxp}{\ensuremath{\mathit{DeltaDnaPaxProt}}}
\newcommand{\deltadmafap}{\ensuremath{\mathit{DeltaDnaMAFAProt}}}
\newcommand{\deltaddeltap}{\ensuremath{\mathit{DeltaDnaDeltaProt}}}

\section{Reactions of the multi-attractor model}
\label{app:A}
The multi-attractor model involves the three species 
$\mafap$, $\deltap$, and $\paxp$ that represent the proteins of the
three genes and it involves ten species that represent the state of the genes:
$\paxd$, $\mafad$, $\deltad$, 
$\paxddeltap$, $\mafadpaxp$, $\mafadmafap$, $\mafaddeltap$,
$\deltadpaxp$, $\deltadmafap$, $\deltaddeltap$.
The chemical reactions are as follows:

\begin{equation}
	\begin{aligned}[c]
     \paxd &\ \stackrel{c_p}{\longrightarrow}\ \paxd + \paxp\\
     \paxp &\ \stackrel{c_d}{\longrightarrow}\ \emptyset \\
     \paxd + \deltap &\ \stackrel{c_b}{\longrightarrow}\ \paxddeltap\\
     \paxddeltap &\ \stackrel{c_u}{\longrightarrow}\ \paxd + \deltap\\
     \mafad &\ \stackrel{c_p}{\longrightarrow}\ \mafad + \mafap\\
     \mafap &\ \stackrel{c_d}{\longrightarrow}\ \emptyset \\
     \mafad + \paxp &\ \stackrel{c_b}{\longrightarrow}\ \mafadpaxp\\
     \mafadpaxp &\ \stackrel{c_u}{\longrightarrow}\ \mafad + \paxp\\
     \mafadpaxp &\ \stackrel{c_p}{\longrightarrow}\ \mafadpaxp + \mafap\\
     \mafad + \mafap &\ \stackrel{c_b}{\longrightarrow}\ \mafadmafap\\
     \mafadmafap &\ \stackrel{c_u}{\longrightarrow}\ \mafad + \mafap\\
     \mafadmafap &\ \stackrel{c_p}{\longrightarrow}\ \mafadmafap + \mafap\\
     \mafad + \deltap &\ \stackrel{c_b}{\longrightarrow}\ \mafaddeltap\\
     \mafaddeltap &\ \stackrel{c_u}{\longrightarrow}\ \mafad + \deltap\\
     \deltad &\ \stackrel{c_p}{\longrightarrow}\ \deltad + \deltap\\
     \deltap &\ \stackrel{c_d}{\longrightarrow}\ \emptyset \\
     \deltad + \paxp &\ \stackrel{c_b}{\longrightarrow}\ \deltadpaxp\\
     \deltadpaxp &\ \stackrel{c_u}{\longrightarrow}\ \deltad + \paxp\\
     \deltadpaxp &\ \stackrel{c_p}{\longrightarrow}\ \deltadpaxp + \deltap\\
     \deltad + \mafap &\ \stackrel{c_b}{\longrightarrow}\ \deltadmafap\\
     \deltadmafap &\ \stackrel{c_u}{\longrightarrow}\ \deltad + \mafap\\
     \deltad + \deltap &\ \stackrel{c_b}{\longrightarrow}\ \deltaddeltap\\
     \deltaddeltap &\ \stackrel{c_u}{\longrightarrow}\ \deltad + \deltap\\
     \deltaddeltap &\ \stackrel{c_p}{\longrightarrow}\ \deltaddeltap + \deltap\\
	\end{aligned}
\end{equation}

\end{appendix}

\end{document}